\documentclass{gtmon_a}
\pdfoutput=1

\proceedingstitle{Proceedings of the School and Conference in Algebraic
Topology (The Vietnam National University, Hanoi, 9-20 August 2004)}
\conferencestart{9 August 2004}
\conferenceend{20 August 2004}
\conferencename{School and Conference in Algebraic Topology}
\conferencelocation{Vietnam National University, Hanoi, Vietnam}

\editor{John Hubbuck}
\givenname{John}
\surname{Hubbuck}

\editor{Nguy\~\ecircumflex{}n H V H\uhorn{}ng}
\givenname{H\uhorn{}ng}
\surname{Nguy\~\ecircumflex{}n}

\editor{Lionel Schwartz}
\givenname{Lionel}
\surname{Schwartz}

\title{Unstable modules over the Steenrod algebra revisited}

\author{Geoffrey M\,L Powell}
\givenname{Geoffrey M L}
\surname{Powell}
\address{Laboratoire Analyse G{\'e}om{\'e}trie et Applications\\
UMR 7539\\\newline
Institut Galil{\'e}e\\
Universit{\'e} Paris 13\\
93430 Villetaneuse\\
France}
\email{powell@math.univ-paris13.fr} 
\urladdr{http://www.math.univ-paris13.fr/~powell/}

\volumenumber{11}
\issuenumber{}
\publicationyear{2007}
\papernumber{12}
\startpage{245}
\endpage{288}

\doi{}
\MR{}
\Zbl{}

\arxivreference{}

\keyword{Steenrod algebra}
\keyword{unstable module}
\keyword{tensor abelian category}
\keyword{functor category}
\subject{primary}{msc2000}{55S10}
\subject{secondary}{msc2000}{18E10}

\received{10 February 2005}
\revised{1 July 2005}
\accepted{24 October 2005}
\published{14 November 2007}
\publishedonline{14 November 2007}
\proposed{}
\seconded{}
\corresponding{}
\editor{}
\version{}

\input xypic
\let\xysavmatrix\xymatrix
\def\xymatrix{\disablesubscriptcorrection\xysavmatrix}
\AtBeginDocument{\let\bar\wbar\let\tilde\wtilde\let\hat\what}

\newcommand{\anarrtilde}{\tilde{\mathcal{A} ^*}}

\renewcommand{\hom}{\mathrm{Hom}}
\newcommand{\calc}{\mathcal{C}}

\newcommand{\calb}{\mathcal{B}}
\newcommand{\cala}{\mathcal{A}}
\newcommand{\calo}{\mathcal{O}}
\newcommand{\cals}{\mathcal{S}}
\newcommand{\calt}{\mathcal{T}}
\newcommand{\cald}{\mathcal{D}}
\newcommand{\bbar}{\wwbar{\calb}}
\newcommand{\atilde}{\wwtilde{\cala}^*}

\newcommand{\spec}{\mathrm{Spec}}

\newcommand{\comod}{\mathit{Comod}}
\newcommand{\comodgr}{{\mathit{Comod}}^{\mathsf{gr}}}
\newcommand{\module}{\mathsf{Mod}}
\newcommand{\vs}{\mathcal{V}_{\field}}
\newcommand{\vsgrad}{{\vs^{\star}}}
\newcommand{\vsbigrad}{{\vs^{\star \star}}}
\newcommand{\vsgr}{\vs^{\zed/2}}
\newcommand{\vsfd}{\vs^{\mathrm{fd}}}
\newcommand{\field}{\mathbb{F}}

\newcommand{\zed}{\mathbb{Z}}

\newcommand{\op}{^\mathrm{op}}

\newcommand{\superalg}{\smash{\mathit{Alg}^{\zed/2}_\field}}
\newcommand{\superbialg}{\smash{\mathit{Bialg}^{\zed/2}_\field}}

\newcommand{\unit}{\mathbf{1}}
\newcommand{\rep}[1][]{\mathrm{Rep}_{#1}}
\newcommand{\hctensor}{\hat{\otimes}}
\newcommand{\comptensor}{\hat{\otimes}}
\newcommand{\fend}{\mathrm{End}_{\f}}
\newcommand{\hf}{\mathrm{Hom}_{\f}}
\newcommand{\f}{\mathcal{F}}

\newcommand{\unst}{\mathcal{U}}
\newcommand{\unstbigr}{{\unst^{\mathrm{bi.gr}}}}
\newcommand{\unstbigrt}{{\unst^{\mathrm{bi.gr},t}}}
\newcommand{\unstmot}{\unst_{\mathcal{M}}}
\newcommand{\hmot}{H_{\mathcal{M}}}
\newcommand{\unsteven}{\unst'}

\makeatletter
\def\cnewtheorem#1[#2]#3{\newtheorem{#1}{#3}[subsection]
\expandafter\let\csname c@#1\endcsname\c@thm}

\def\dnewtheorem#1[#2]#3{\newtheorem{#1}{#3}
\expandafter\let\csname c@#1\endcsname\c@THM}

\newtheorem{thm}{Theorem}[subsection]
\cnewtheorem{prop}[thm]{Proposition}
\cnewtheorem{cor}[thm]{Corollary}
\cnewtheorem{lem}[thm]{Lemma}
\newtheorem{THM}{Theorem}
\dnewtheorem{COR}[THM]{Corollary}
\dnewtheorem{PROP}[THM]{Proposition}

\theoremstyle{definition}
\cnewtheorem{defn}[thm]{Definition}
\cnewtheorem{exam}[thm]{Example}

\theoremstyle{remark}
\cnewtheorem{rem}[thm]{Remark}
\cnewtheorem{nota}[thm]{Notation}
\cnewtheorem{hyp}[thm]{Hypothesis}

\makeatother  

\makeautorefname{hyp}{Hypothesis}

\begin{document}

\begin{abstract}
A new and natural description of the category of unstable modules over
the Steenrod algebra as a category of comodules over a bialgebra is 
given; the theory extends and unifies the work  of Carlsson, Kuhn, 
Lannes, Miller, Schwartz, Zarati and others.  Related categories of
comodules are studied, which shed light upon the structure of the 
category of unstable modules at odd primes. In particular, a category 
of bigraded unstable modules is introduced; this is related to the study 
of modules over the motivic Steenrod algebra.
\end{abstract}

\maketitle

\section{Introduction}
\label{section-introduction}

The Steenrod algebra $\cala$ over a prime field $\field$ of
characteristic $p$ is a fundamental mathematical object; it is defined
in algebraic topology to be the algebra of stable cohomology
operations for singular cohomology with coefficients in $\field$. The
algebra $\cala$ is graded and acts on the cohomology ring of a space;
the underlying graded $\cala$--module is {\em unstable}, a condition
which is usually defined in terms of the Steenrod reduced power
operations and the Bockstein operator. 

This  paper shows that  the well-known description
of the dual of the Steenrod algebra, which is due to Milnor, has an 
extension which leads to an entirely algebraic description of the category 
$\unst$  of
unstable modules over the Steenrod algebra as a category of comodules 
(defined with respect to a completed tensor product) 
over a  bialgebra, without imposing an external instability condition. The 
existence of such a description is
implied by the general theory of tensor abelian categories at the prime 
two; in the odd prime  case, the method requires the  super-algebra 
setting, namely 
using $\zed/2$--gradings to introduce the necessary sign conventions for
commutativity. 

\begin{THM}
The category of unstable modules over the Steenrod algebra at  an odd
prime is equivalent to the category of  right comodules (with respect
to a completed tensor product) over a  $\zed/2$--graded bialgebra
$\atilde$.
\end{THM}

At an odd prime,  the theory introduces an auxiliary  bialgebra $\calb$,  
which  can be defined as the bialgebra of endomorphisms of the additive 
group (in the super-algebra setting). The
category $\unst (\calb)$ of graded comodules over this coalgebra is
equipped with an exact forgetful functor $\unst (\calb) \rightarrow
\unst$ to the category of unstable modules,  which is induced by a
surjective morphism of bialgebras $\calb \twoheadrightarrow \atilde$. 
Moreover, $\unst(\calb)$ is
equivalent to a representation category which generalizes the
description given by Kuhn of the category of unstable modules over the
 field $\field_2$. This part of the theory extends and unifies existing 
approaches to the category of unstable modules. 

The surjection  $\calb \twoheadrightarrow \atilde$ factorizes across a 
bialgebra $\bbar$; the category
$\unstbigr$ of comodules over this bialgebra is a bigraded analogue of
the category of unstable modules. This category has not hitherto been 
studied.  The structure of $\unstbigr$ is of
interest since it sheds light on and provides new approaches to the 
structure of the classical category of unstable modules at odd primes and 
also in connection with unstable modules defined over the motivic Steenrod 
algebra (at all
primes). At an odd prime,  these tensor abelian categories are related by 
exact functors which are induced by corestriction
\[
\unst (\calb) 
\stackrel{\Psi}{\rightarrow} 
\unstbigr
\stackrel{\Theta}{\rightarrow}
\unst.
\]
The functorial point of view  on the category of unstable modules,
developed by Kuhn \cite{k} from the work of Henn, Lannes and Schwartz
\cite{hls} 
arises naturally in the theory via the free commutative algebra
functor; forgetting the algebra structure, this can be considered as
an object of the category $\f$ of functors  from finite-dimensional
$\field$--vector spaces to $\field$--vector spaces. Key examples of
objects of $\f$ are given by the exterior power functors $\Lambda ^a$
and the divided power functors $\Gamma ^b$.

The one-sided Morita equivalence theory of \cite{k} gives rise to  
functors $r_\calb \co \f \rightarrow \unst (\calb)$ and  
$ r ' \co \f \rightarrow \unsteven$, 
where $\unsteven \subset \unst$ denotes the full 
subcategory of unstable modules concentrated in even degree. These give a 
very natural description of the projective generators of the category 
$\unstbigr$.

 The results on the projective generators of $\unstbigr$  lead,  at odd
 primes, to a
 new and natural analysis of the free unstable  modules $F(n)$. The
 following result provides an  odd-primary analogue of the
 well-known analysis of the structure of the free unstable  modules
 at the prime two, where $\calo$ denotes the forgetful functor 
$\unsteven \rightarrow \unst$. 

\begin{PROP}
For $p$ an odd prime 
and $n$ a positive integer,  $F(n)$ has a finite filtration with  
associated graded 
\[
\bigoplus_{\substack{a+2b=n\\a \geq 0}}
\calo r ' (\Lambda ^a \otimes \Gamma ^b) 
\oplus
\bigoplus
_{\substack{a+2b=n\\
a \geq 1} }
\Sigma \calo r' ( \Lambda ^{a-1} \otimes \Gamma ^b ).  
\]
 \end{PROP}

There is a corresponding analysis of the injective cogenerators of
$\unst$, which is presented in \fullref{subsect-cogen-inj}. These
results are applied in  \fullref{sect:new-proofs}  to give new proofs of 
the fundamental results on nil-localization in the odd characteristic
case, using the results on $\unst (\calb)$ provided by the
representation-theoretic framework. 

The final section of the paper indicates a modification of the theory
at the prime two which introduces an abelian category $\unstmot$ which
is related to the study of unstable  modules over the motivic Steenrod
algebra. This allows the comparison of the category of unstable
modules and a suitable category of bigraded unstable modules at the
prime two. 

%
\subsection{Related results}

There are related results which occur in the literature; the results on 
unstable modules at the prime two are implicit in the work of Kuhn 
\cite{k}. The bialgebra which is used to define the category of unstable 
modules occurs in Bisson--Joyal \cite[Section 4]{bj}, where it is 
termed the extended Milnor--Hopf algebra. 
The addendum \cite[page~260]{bj} indicates 
the fact that the category of unstable modules corresponds to the category 
of comodules over this bialgebra. A second reference for related material 
in a similar context is Smirnov \cite[page~116, Chapter 5]{sm}. This
reference also provides a related statement for the odd prime case.

\subsection{Organization of the paper}

The first part of the paper is devoted to the introduction of the 
categories of (generalized) unstable modules which are of interest here. 
\fullref{section-tensor-abelian} provides a survey of the theory of tensor 
abelian categories which motivates the constructions of the paper. 
\fullref{section-universal-bialgebra} constructs the bialgebra $\calb$, by 
considering the endomorphisms of the additive group; this is a 
generalization of Milnor's approach to calculating the dual Steenrod 
algebra. The other bialgebras used in the paper are constructed as 
quotients of  $\calb$. \fullref{section-comodules-unstable} defines the 
categories of generalized unstable modules as categories of  comodules; 
the simple objects of these categories are considered briefly, together 
with associated suspension functors.  

The second part of the paper establishes the connection with the 
functorial point of view. \fullref{section-functors-bialgebras} reviews 
the category $\f$ of functors, the notion of an exponential functor and 
then establishes the relation with the bialgebra $\calb$ defined in the 
first part. \fullref{sect-rep-cat} reviews and extends the results of Kuhn 
on representation categories, motivated by the considerations of this 
paper. 

The third part  is devoted to an analysis of the projective and injective 
objects of the category $\unstbigr$ at an odd prime. 
\fullref{section-proj-inj} considers the standard projective generators 
and the standard injective cogenerators of $\unstbigr$; the results are 
new and are applied to obtain a new analysis of the structure of the 
projective and injective objects of the category $\unst$ of unstable 
modules. \fullref{sect:new-proofs}  uses the theory to provide a new proof 
of the injectivity of $H^* (BV)$ in the odd prime case, more in the spirit 
of \cite{k}.

The fourth part of the paper corresponds to \fullref{sect-motivic}; this 
indicates how bigraded unstable modules appear at the prime two. This 
theory is related to the study of modules over the motivic Steenrod 
algebra. This material will be returned to in greater detail elsewhere. 

The appendix reviews certain results on comodules which are required in 
the paper. 

\part{Basic structure}

%
%
\section{Tensor abelian categories}
\label{section-tensor-abelian}

The category of unstable modules over the Steenrod algebra  is a
tensor abelian category to which   Tannakian theory can be
applied. The relevant theory of tensor abelian categories is reviewed
in this section;  for further details, the reader is referred to 
Deligne--Milne \cite{dm}.

\subsection{General theory}

A tensor category is a category $\calc$ which is equipped with a
symmetric monoidal structure $(\calc, \otimes ,\unit )$, where $\unit$
denotes the unit.  For a field $\field$,  let $\vs$ denote the category of 
$\field$--vector spaces, equipped with the usual abelian tensor structure.


\begin{defn}
\ 
\begin{enumerate}
\item
An  abelian tensor category is a tensor category $(\calc ,
\otimes , \unit)$ for which the category $\calc$ is abelian and the
functor $\otimes $ is biadditive. 
\item
For $(\calc, \otimes , \unit)$ an  abelian tensor category  such
that $\mathrm{End} (\unit) = \field$ is a field,  a  fibre functor
 is a faithful, exact, $\field$--linear tensor functor $\calc
\rightarrow \vs$. 
\end{enumerate}
\end{defn}

An affine monoid scheme over $\field$  is a scheme of the form 
$M := \spec (B)$,
where $B$ is a bialgebra over $\field$  (not necessarily equipped with a
conjugation), for which the underlying $\field$--algebra is
commutative\footnote{This usage of the term bialgebra conflicts with the
terminology of \cite{dm}.}.  An affine group scheme is an affine monoid  
scheme of the form $G : = \spec (H)$,
where $H$ is a Hopf algebra over $\field$  (a bialgebra with
conjugation).

\begin{nota}
 For $M = \spec (B)$ an affine monoid scheme over $\field$, let  
$\rep[\field] (M)$ (the category  of representations of $M$) denote the 
category of right $B$--comodules. 
\end{nota} 

\begin{prop}
\label{prop:rep-ta}
For $M$  an affine monoid scheme, the category $\rep[\field](M)$ is a
tensor abelian category, equipped with a canonical fibre functor.
\end{prop}

\begin{exam}
Let $\field $ be a field.
\begin{enumerate}
\item 
The category of comodules over the  multiplicative group $\mathbb{G}_m 
\cong \spec (\field[x^{\pm 1}])$ is equivalent to  the category of graded 
vector spaces. 
\item
The category of comodules over  $G := \spec (\field [t]/
(t^2 -1))$ is equivalent to the category of  $\zed/2$--graded vector 
spaces. 
\item
\cite[Example 1.25]{dm} Let $\vsgr$ denote the category of
$\zed/2$--graded vector spaces, equipped with the symmetric monoidal
tensor structure involving the Koszul sign convention. Then $\vsgr$ is
a tensor abelian category but it is not of the form $\rep[\field]
(M)$, 
for any affine monoid scheme $M$.
\end{enumerate}
All of the above examples are {\em rigid} tensor categories, which means 
that duality behaves well. In the first two cases, this follows from the 
fact that the categories are representations of affine group schemes 
rather than just affine monoid schemes.
\end{exam}
  
The following is the part of the theory  of  neutral Tannakian
categories which is relevant to this paper.

\begin{thm}
[Deligne and Milne {{\cite[Proposition 2.14]{dm}}}]
\label{thm:Tann}
Let $(\calc , \otimes , \unit)$ be a $\field$--linear abelian tensor
category such that $\mathrm{End} (\unit) = \field$ and let 
$\omega \co \calc \rightarrow \vsfd$ be a fibre functor to the category of 
finite-dimensional vector spaces. There exists an affine monoid scheme $M 
= \spec (B)$ over $\field$
such that $(\calc , \otimes, \unit)$ is equivalent to the category
 $\rep[\field](M)$  equipped with the canonical abelian tensor category
structure and fibre functor.
\end{thm}

\begin{rem}
The above theory can be generalized to the $\zed/2$--graded situation
with the Koszul sign convention. In this context, the fibre functor
is a tensor functor 
$
\omega _\calc \co  \calc \rightarrow \vsgr.
$

This generalization is essential for the topological considerations,
where the algebras to be considered are naturally $\zed/2$--graded and
are commutative with respect to the Koszul sign convention.
\end{rem}

\begin{exam}
Let $\field$ be the prime field of characteristic $p$ and let $\unst$
denote the category of unstable modules over the $\field$--Steenrod
algebra. 
\begin{enumerate}
\item
For $p=2$, the category $\unst$ is an $\field_2$--linear
abelian tensor category, equipped with a fibre functor to the category
of $\field_2$--vector spaces. 
\item
For $p>2$, the category $\unst$ is an $\field_p$--linear abelian tensor
category, equipped with a fibre functor to the category $\vsgr $ of
$\zed/2$--graded $\field_p$--vector spaces, equipped with the tensor
structure with the Koszul sign convention.
\end{enumerate}
\end{exam}

The theory of Tannakian categories is  developed in terms of
finite-dimensional representations. For this reason, Tannakian theory
does not apply directly in considering unstable modules, but requires 
modification using completed
tensor products in the comodule structures. Modulo this addendum,
Tannakian theory implies 
 that the category of unstable modules (for $p=2$ and for $p$
odd) has a description as a category of comodules.

\section{Endomorphisms of the additive group}
\label{section-universal-bialgebra}

This section constructs the bialgebra $\calb$, which is the
fundamental mathematical object of this paper, as  the endomorphism
bialgebra of the additive group in the category of $\zed/2$--graded
algebras. Throughout the section, $\field$ denotes a prime field of
odd characteristic.

The $\field$--algebra structures considered in this
paper are graded commutative, with  respect to the Koszul sign
convention. The foundations are developed in the super-algebra
context, namely using 
$\zed/2$--gradings, to avoid  imposing an {\em a priori}  
$\zed$--grading.  The fact that the foundations of algebraic geometry can 
be generalized to the  super-algebra context is well known (see Deligne
\cite[Section 0.3]{d2}). 
\subsection{Super algebras}
Let $\vsgr$ denote the category of $\zed/2$--vector spaces, equipped
with the symmetric monoidal structure provided by the graded tensor
product with the  Koszul sign convention \cite[Example
1.25]{dm}. Let $\superalg$ denote the category of unital commutative 
monoids
in $\vsgr$; this is the category of unital $\zed/2$--graded
$\field$--algebras which are graded commutative, with respect to the
Koszul sign convention. The tensor product of $\vsgr$ induces the
coproduct in the category $\superalg$; in particular, the tensor
product of two objects of $\superalg$ is in $\superalg$. 

The category $\superbialg$ of $\zed/2$--graded bialgebras is the
category of comonoid objects in $\superalg$. Namely an object $B \in
\superbialg$ is a $\zed/2$--graded algebra $B \in \superalg$ which is
equipped with a morphism\footnote{Here, as elsewhere, the graded
tensor product is denoted simply $\otimes$.}
$B \rightarrow B \otimes B$ in $\superalg$,
which is coassociative and counital with respect to the counit morphism $B 
\rightarrow \field$ in
$\superalg$. 

%
\subsection[The bialgebra]{The bialgebra $\calb$}

\begin{nota}
Let $H$ denote the free $\zed/2$--graded algebra on the $\zed/2$--graded
vector space $\langle x, y \rangle$, where $y$ has degree $1$ and $x$
has degree $0$.
\end{nota}

The algebra $H$ has the structure of a $\zed/2$--graded Hopf algebra,
with underlying algebra $\Lambda (y) \otimes \field[x]$, which is
primitively-generated. The Hopf algebra $H$ represents the additive
group in the context of $\zed/2$--graded algebras. 

\begin{rem}
The algebra $H$ has a Hausdorff  filtration $(I^n H)$ given  by powers of 
the
augmentation ideal, hence it is possible to form half-completed tensor
products $H \hctensor V:= \lim_\leftarrow (H/ I^n H) \otimes V $, for
any vector space $V$.

The half-completed tensor product  is sufficient
to be able to define a general notion of completed comodule
structure over a bialgebra $B$, where the structure morphism is given
by a morphism $\psi \co H \rightarrow H \hctensor B$. 

In the applications, all objects will have a $\zed$--grading and all
completed tensor products can be understood in the usual graded
context. For this reason, the details concerning the usage of the
half-completed tensor product are left to be supplied by the
interested reader. 
\end{rem} 

All comodule structures in this section are understood to be defined
with respect to a half completed tensor product.

\begin{defn}
For $B\in \superbialg$ a $\zed/2$--graded bialgebra, a multiplicative
right $B$--comodule algebra structure  on $H$ is a morphism of 
super-algebras $\psi \co H \rightarrow H \hat{\otimes} B$ which induces a 
comodule structure on $H$. 
\end{defn}

\begin{prop}
\label{prop:calb-construct}
\ 
\begin{enumerate}
\item
There exists a $\zed/2$--graded bialgebra $\calb$ with underlying
$\field$--algebra the free super-algebra
$$
\calb \cong \Lambda (w, \tau_i|i \geq 0) \otimes \field [u, \xi_j |j \geq 
0]
$$
and with coproduct
\begin{eqnarray*}
 u
 &\mapsto &
u \otimes u + w \otimes \tau_0\\
w
 &\mapsto&
u \otimes w + w \otimes \xi_0\\
 \tau_i
 &\mapsto&
\sum_{s=0} ^i  \xi_{i-s} ^{p^s} \otimes \tau_s + \tau_i \otimes u\\
 \xi_j 
&\mapsto&
\sum _{s =0} ^j 
   \xi_{j-s} ^{p^s} \otimes \xi_s + \tau _j \otimes w.
\end{eqnarray*}
\item
The underlying super--algebra $H$ admits a $\calb$--comodule
structure $\psi \co H \rightarrow H \hctensor \calb$ such that 
\begin{enumerate}
\item
$\psi$ is a morphism of super-algebras, determined by 
$$y
\mapsto
y \otimes u + \sum _{i \geq 0} x^{p^i} \otimes \tau_i
\quad\text{and}\quad
x 
\mapsto
y \otimes w + \sum _{j _\geq 0}  x^{p^j} \otimes \xi_j;$$
\item
the coproduct $H \rightarrow H \otimes H$ is a morphism of
$\calb$--comodules, where $H \otimes H$ is given the tensor product
 $\calb$--comodule structure.  
\end{enumerate}
\end{enumerate}
\end{prop}

\begin{proof} (Indications) \qua The construction of the bialgebra, its 
coproduct
and the comodule structure is a straightforward generalization of
Milnor's method \cite{mi} for  calculating  the dual of the Steenrod 
algebra.
\end{proof}

\begin{rem}
\ 
\begin{enumerate}
\item
The bialgebra $\calb$ can be interpreted as the endomorphism bialgebra
of the additive algebraic group; this  implies that it satisfies a 
universal
property, the formulation of which is left to the reader.
\item
The case of the prime field of characteristic two is similar, but more
elementary, since the algebras considered are commutative in the
ungraded sense and the respective Hopf algebra $H$  is the polynomial
algebra on a single generator. The universal bialgebra $\calb _2$ is 
$\field_2
[\xi_j | j \geq 0]$, equipped with the coproduct
$$\Delta \xi_j =
\sum_{s=0}^j \xi_{j-s}^{2^s} \otimes \xi_s.$$
This is the extended Milnor--Hopf algebra of Bisson--Joyal \cite[Section 4]{bj}.
\item
It is well-known that the Steenrod algebra, for $p=2$, can be regarded
as automorphisms of the additive formal group  and it is folklore that
this can be extended, for $p$ odd,  by considering the super-algebra
setting. The above generalizes this point of view,  by considering the 
full endomorphism ring in the $\zed/2$--graded setting.
\end{enumerate}
\end{rem}

\subsection
[Bialgebras derived from the universal bialgebra B]{Bialgebras derived 
from the universal bialgebra $\calb$ }There are quotient bialgebras of
$\calb$ which are of importance in considering the category of
unstable modules over the Steenrod algebra. The ideals $\langle w \rangle$ 
and $\langle w, \xi_0 - u ^2 \rangle $ are Hopf ideals in $\calb$, which 
allows the following definition.

\begin{defn}
\quad
\begin{enumerate}
\item
Let $\bbar$ denote the quotient  $\calb / (w)$ in
$\superbialg$.  
\item
Let $\atilde$ denote the quotient a $\calb / (w, \xi_0 - u^2 )$ in
$\superbialg$.\footnote{The notation reflects the relation with the dual  
of the Steenrod algebra.} 
\item
Let $\calb'' \in \superbialg$ denote the sub-bialgebra of $\atilde$ which
is generated by the elements $\xi_j$ together with $u$. 
\end{enumerate}
\end{defn}

\begin{rem}
\ 
\begin{enumerate}
\item
The bialgebras $\calb, \bbar, \atilde \in \superbialg$ do not have 
$\zed/2$--graded Hopf
algebra structures. For example, $\field [\xi_0]$ is a sub-bialgebra
of $\bbar$ which does not have the structure of a Hopf algebra, since 
$\xi_0$ is grouplike and not invertible. 
\item
The dual of the Steenrod algebra, $\cala ^*$,  is obtained as the quotient 
of the bialgebra $\atilde$ by the Hopf ideal generated by $u 
-1$.\footnote{It is the desire for compatibility with the usual notation 
for the dual of the Steenrod algebra which imposes the cumbersome notation 
for $\atilde$.}  
\item
The bialgebra $\atilde$ is related to the bigraded algebra $J ^* _*$ which 
was constructed by Miller from the Brown--Gitler modules (cf Schwartz
\cite[Theorem 2.4.8]{s}).
\end{enumerate}
\end{rem}

\begin{lem}
\label{lem:comm-bialgebras}
There is a commutative diagram of $\zed/2$--graded bialgebras:
\[
\xymatrix{
\calb \ar[r]
&
\bbar 
\ar[r]
\ar[d]
&
\atilde
\ar[d]
\\
&
\field [\xi_0, u ] 
\ar[r]_{\xi_0 \mapsto u^2} 
&
\field [u] 
}
\]
where the bialgebras $\field [\xi_0, u], \field [u]$  are generated by
group-like elements. Moreover,  the vertical morphisms of the diagram are 
split
surjections of bialgebras.  
\end{lem}

\subsection[Introducing z gradings]{Introducing $\zed$--gradings}

The category of non-negatively graded $\field$--vector spaces is 
equivalent to the category of right comodules over $\field[u]$. The 
following results show that there are natural gradings on the bialgebras 
under consideration, which induce the $\zed/2$--gradings. 

\begin{lem}
\label{lem-left-finite-type}
For   $B$ one of the bialgebras $\calb, \bbar, \calb'', \atilde,$
\begin{enumerate}
\item $B$ admits a morphism of bialgebras $B \rightarrow \field [u]$,
  which factors across  $\atilde
\rightarrow \field [u]$;
\item
 $B$ is naturally bigraded with
respect to the corestricted  left and right  $\field [u]$--comodule
structures;
\item
$B$ is of finite
type with respect to the grading induced by the left comodule structure.
\end{enumerate}
\end{lem}

\begin{proof}
The morphism $B \rightarrow \field [u]$ is provided by 
\fullref{lem:comm-bialgebras}; the remainder of the Lemma is
straightforward. 
\end{proof}

\begin{defn}
Let $B$ be a bialgebra as above, equipped with the bigrading induced by 
the $\field[u]$--comodule structures; the total degree of a bihomogeneous 
element of bidegree $(m,n)$ is the integer $m-n$. 
\end{defn}

The following result implies that all the $\zed/2$--gradings which are
considered are the reduction of a natural $\zed$--grading on the
algebra $\calb$. There are similar results for the other bialgebras.

\begin{prop}
\label{prop:zed-compat}
The $\zed/2$--grading of $\calb$ is the mod $2$ reduction of the
$\zed$--grading of $\calb$ given by the total degree. 
\end{prop}

\begin{rem}
The usual grading of the dual Steenrod algebra $\cala^*$ can be
recovered from the total degree, using the total degree defined on
$\atilde$ together with the observation that the total degree of $u$
is zero. 
\end{rem}

\section{Categories of comodules related to unstable modules}
\label{section-comodules-unstable}

This section defines the categories of comodules which are of interest
in this paper, in relation to the category of unstable modules over
the Steenrod algebra.

Throughout this section, let $\field$ be the prime field of
characteristic $p$, where $p >2$; the underlying category of vector
spaces  is taken to be the category of
$\zed/2$--graded vector spaces, equipped with the Koszul--sign tensor
 structure.

The results presented in this section have analogues for the case $p=2$.  

%
\subsection{Categories of graded comodules}
\label{subsect:graded-comod}

\begin{defn}
\label{def-graded-comod}
Let $B \in \superbialg$ be a bialgebra which is equipped with a
 morphism of bialgebras $B \rightarrow \field[u]$ and hence has a left and 
right grading.
 Suppose that $B$ is of finite type with respect to the left grading. 

 Let $\comodgr _ B$ be the category of graded, 
 $B$--comodules, in the following sense. An object of $\comodgr _B$ is
 a non-negatively graded vector space $M$, with grading defined by the
 comodule structure $M \rightarrow M \otimes \field [u]$, together
 with the comodule structure morphism $M \rightarrow M \comptensor B$
 which satisfies the condition that the diagram 
\[
\xymatrix{
M \ar[r] 
\ar[d]
&
M \comptensor B 
\ar[d]
\\
M \otimes \field [u] 
\ar[r]
&
M\comptensor \field [u]
}
\]
commutes. In particular, the morphism $M \rightarrow M \comptensor B$
is a morphism of graded vector spaces, where the grading on the right
is induced by the right grading of $B$.
\end{defn}

\begin{prop} 
Let $B$ be a bialgebra which is equipped with a morphism of bialgebras
$B \rightarrow \field[u]$. Suppose that $B$ is of finite type with
respect to the grading induced by the associated left
$\field[u]$--comodule structure, then the category $\comodgr _B$ is an
abelian tensor category. 
\end{prop}

The bialgebras $\calb , \bbar , \calb'', \atilde $ satisfy the
hypotheses of \fullref{def-graded-comod}, by 
\fullref{lem-left-finite-type}; hence the above definition can be applied.

\begin{defn}
Define the following tensor abelian categories:
\begin{align*}
\unst (\calb )&:= \comodgr _\calb &
\unstbigr &:= \comodgr _{\bbar}\\
\unst''  &:= \comodgr _{\calb''} &
\unst &:= \comodgr _{ \anarrtilde}\\
\unsteven &:= \comodgr _{\field [\xi_j |j  \geq 0]} 
\end{align*}
\end{defn}

The notation $\unst$ does not conflict with the usual usage, by the
following result.

\begin{thm}
\label{thm:unst-comod}
The category $\comodgr _{ \anarrtilde}$ is equivalent to the category of 
unstable modules
over the mod--$p$ Steenrod algebra $\cala$.
\end{thm}

\begin{proof}
The category of unstable modules over the Steenrod algebra is
usually defined as the category of graded modules over the  Steenrod 
algebra, $\cala$, subject to an instability condition in terms of
the operation of elements derived from the dual basis to the basis of
monomials in the elements $\xi_j, \tau_i$. It is elementary to show
that this condition is equivalent to a condition on the adjoint coaction
involving the terms $\xi_1 ^t$. It is a straightforward exercise to
show that this implies that a graded module over the Steenrod algebra
is unstable if and only if the adjoint coaction extends to a graded
right  $\atilde$--comodule structure. 
\end{proof}

The objects of $\unstbigr$ are naturally bigraded (which justifies the
notation), by the following result, in which  $\vsbigrad$ denotes the 
category of $\zed \times \zed$--graded vector spaces.

\begin{prop}
The morphism $\bbar \rightarrow \field[\xi_0, u]$ induces an exact functor 
of tensor abelian categories
$
\unstbigr
\rightarrow 
\vsbigrad.
$
\end{prop}

\begin{nota}
The bidegrees will be written as $(a,b)$, where $a$ denotes the 
$u$--degree and $b$ denotes the $\xi_0$--degree.
\end{nota}

\begin{prop}
There is a diagram of exact functors between abelian categories:
$$
\xymatrix{
&
\unsteven
\ar[r]
\ar@{_(->}[d]_{\hat{\calo}}
&
\unst'' 
\ar@{_(->}[d]
\\
\unst ( \calb) 
\ar[r]_{\Psi}
&
\unstbigr 
\ar[r]_{\Theta}
&
\unst
}
$$
in which the horizontal morphisms indicate forgetful functors. The 
embeddings $\unsteven \hookrightarrow \unstbigr$ and   $\unst'' 
\hookrightarrow \unst$ are 
fully faithful  and admit retractions $\unstbigr \rightarrow \unsteven$ 
and $\unst \rightarrow \unst''$ respectively. 
\end{prop}

\begin{proof}The exact functors are induced by the corestriction functors 
which
are associated to the canonical morphisms of the respective
bialgebras (cf \fullref{lem:comm-bialgebras}). The retractions are given 
by the respective right adjoint functors.
\end{proof}

\begin{rem}
\ 
\begin{enumerate}
\item
The category $\unst''$ is related to  the full sub-category $\unsteven$
of $\unst$ which identifies with the objects which are concentrated in even
degree, which has been used in the work of Lannes and Zarati on the 
category $\unst$. In particular, there is an adjunction, 
$
\mathcal{O} \co \unsteven \rightleftarrows \unst \thinspace{:} \tilde{\mathcal{O}},
$ 
where $\mathcal{O}$ denotes the forgetful functor and 
$\tilde{\mathcal{O}}$ 
its right adjoint. The category $\unst''$ splits as a product of two 
copies of $\unsteven$,
corresponding respectively to elements in even (resp. odd)  degrees.
\item
The category $\unst (\calb)$ allows for the action of the element $w$ of 
negative total degree. 
\item
The category $\unstbigr$ sheds light on the structure of $\unst$; it
is also of interest in  studying  unstable modules over the motivic
Steenrod algebra, since motivic cohomology is naturally bigraded by
the topological degree and the twist (or weight). 
\end{enumerate}
\end{rem}

\begin{exam}
The object $H: = \Lambda (y ) \otimes \field [x]$ has the structure of an 
object in $\unst (\calb ) $ and hence of an object in
$\unstbigr$, by corestriction. Thus, for each
non-negative integer $m$, $H ^{\otimes m}$ defines an object of $\unst
(\calb)$  and therefore  of $\unstbigr$, by corestriction.
\end{exam}

\subsection{Simple comodules}

The simple objects of the abelian categories $\unst$, $\unstbigr$, 
$\unsteven$,
$\unst ''$ are understood via the following result.

\begin{prop}
A simple object $M$  in one of the abelian categories $\unst$,
$\unstbigr$, $\unsteven$, $\unst ''$  has  underlying vector spaces of 
total dimension one, with  comodule structure which corresponds to the 
grading. 
\end{prop}
 
\begin{exam}
The bigraded vector space $\field (a, b)$ of total dimension one,
concentrated in bidegree $(a, b)$ where $a, b$ are non-negative
integers has the natural structure of an object of $\unstbigr$.
\end{exam}

\begin{rem}
The bigraded vector space $\field(a,b)$ does not in general have the
structure of an object of $\unst (\calb)$. The simple objects of the 
category $\unst (\calb)$ are constructed
from the simple objects of the category $\comodgr_{\calc}$, where
$\calc$ denotes the quotient bialgebra of $\calb$ with underlying
algebra $\field [\xi_0, u] \otimes \Lambda
(\tau_0, w)\in \superbialg$.
\end{rem}

\subsection{Suspension functors}

The simple objects of the categories $\unstbigr, \unst , \unsteven, \unst 
''$ define suspension functors by forming the tensor product. 

\begin{defn}
Let $a,b, n$ be non-negative integers. 
\begin{enumerate}
\item 
The suspension functor $\Sigma^{(a,b)}$ of bidegree $(a,b)$ is the functor 
defined by the tensor product 
 $\field (a, b) \otimes - \co \unstbigr\rightarrow \unstbigr$. 
\item
The suspension functor $\Sigma ^n \co \unst \rightarrow \unst $ is the 
functor defined by the tensor product $\field (n ) \otimes - \co \unst 
\rightarrow \unst$. 
\end{enumerate}
\end{defn}

The forgetful functor $\Theta \co \unstbigr \rightarrow \unst$ sends $\field
(a,b)$ to the vector space $\field(a+2b)$ of total  dimension one, 
concentrated in degree
$a+2b$. The following lemma is clear.

\begin{lem}
Let $a,b$ be non-negative integers, then there is a commutative
diagram
$$
\xymatrix{
\unstbigr
\ar[r]^{\Sigma^{(a,b)}}
\ar[d]_ \Theta
&
\unstbigr 
\ar[d]^\Theta
\\
\unst 
\ar[r]_{\Sigma ^{a+2b}}
&
\unst.
}
$$
\end{lem}

\part{The functorial viewpoint}
\section{Functors and bialgebras}
\label{section-functors-bialgebras}

The bialgebra $\calb$ is the endomorphism bialgebra of the free
graded-commutative algebra on the $\zed/2$--graded vector space
$\langle x, y \rangle$; forgetting the algebra structure, the free
graded-commutative algebra can be
considered as a functor from $\field$--vector spaces to $\field$--vector
spaces. The relation between the category of unstable modules and the
category of functors between $\field$--vector spaces follows from the
comparison between the bialgebra of endomorphisms of the additive
group and the bialgebra of endomorphisms of the free
graded-commutative algebra functor, which is  presented in 
\fullref{thm:Bend}.

\subsection{Functors}

The category $\f$ of functors from finite-dimensional $\field$--vector
spaces to $\field$--vector spaces is an abelian tensor category with
enough projectives and enough injectives, given respectively by
Yoneda's lemma and its dual.  An object of $\f$ (a functor) is {\em
finite} if it has a finite composition series and is {\em
polynomial} if it is polynomial in the sense of Eilenberg--MacLane
(see Kuhn \cite{k}); these two conditions are equivalent for functors which
take finite-dimensional values. A functor is {\em analytic} if it is
the colimit of its finite subobjects; the full subcategory of analytic 
functors is denoted $\f_\omega$. 

 The divided power functors, $\Gamma ^n$, defined by $V \mapsto (V 
^{\otimes n } ) ^{\mathfrak{S}_n}$,  the symmetric power functors, $S^n$, 
and the exterior power functors, $\Lambda ^n$, are finite functors which 
are of fundamental importance to the theory. The embedding theorem of 
\cite{k} is interpreted as follows.

\begin{thm}[Kuhn \cite{k}]
\label{thm-emb-thm}
For $F$ a finite functor in $\f$, there exists a finite set of 
non-negative integers $\{n_ i \}$ and a surjection $\bigoplus _i \Gamma 
^{n_i} \twoheadrightarrow F$.
\end{thm}

An important observation which is used in the proof is the following: 
\begin{lem}
\label{lem-emb-tensor}
Suppose that $\Gamma ^a \twoheadrightarrow F_1 $, $\Gamma ^b 
\twoheadrightarrow F_2$ are surjections in $\f$, then there exists  a 
surjection $\Gamma ^N \twoheadrightarrow F_1 \otimes F_2$ in $\f$, for 
some  non-negative integer $N$.
\end{lem}  

An extensive array of calculations of Ext groups in the category $\f$ have 
been performed (see Franjou--Friedlander--Scorichenko--Suslin
\cite{ffss}, for example). For the purposes of this 
paper, the following elementary calculation is important. 

\begin{lem}
For $a,b$ non-negative integers,
\begin{enumerate}
\item
$\hf (\Gamma ^b, \Lambda ^a) = 
\begin{cases}
  \field & b = \sum_{i=1} ^a p^{n_i}, n_i < n_{i+1} \\
0 & \text{otherwise;}
\end{cases}
$
\item
for $p >2$, 
$\hf ( \Lambda ^a, \Gamma ^b) = 
\begin{cases}
  \field & a=b,  a \in \{0, 1 \} \\
0 & \text{otherwise.}
\end{cases}$
\end{enumerate}
\end{lem}

This result, together with the calculation of $\hf (\Gamma ^* , \Gamma^*)$ 
(see \cite{k}),  is used in conjunction with the exponential property of 
the functors $\Gamma ^*, \Lambda ^*$.

\subsection{Generalities on exponential functors}

A functor $E \in  \f$ is  exponential if there exists a binatural
isomorphism $E (V\oplus W) \cong E(V) \otimes E(W)$, for $V, W \in
\vsfd$. A graded functor is exponential if it satisfies this property
with respect to the graded tensor product\footnote{For the remainder
of this section, `exponential functor' is used to indicate either the
graded or the ungraded version - the context should make the meaning
clear.}.

\begin{exam}
The injective functor
$I_V$, for $V \in \vsfd$, defined  by 
$I_V \co W \mapsto \field ^{\hom (W, V)}$  is an exponential functor.
\end{exam}

\begin{rem}
\ 
\begin{enumerate}
\item
The tensor product of two (ungraded) exponential
functors is  an exponential functor; in the graded case, one obtains a 
bigraded exponential functor. 
\item
The structure of an exponential functor $E$ induces a canonical
product $\mu \co E \otimes E \rightarrow E$ and coproduct 
$\Delta \co E \rightarrow E\otimes E$.
\end{enumerate}
\end{rem}

\begin{exam}
For $\field$ a prime field of odd characteristic,  the graded functor
$\Lambda^*$ is a graded exponential functor; the associated product
and coproduct are graded commutative, since  Koszul signs
intervene. More generally, for the topological application of this
paper, it is necessary to consider the  bigraded exponential functor
which corresponds to $\Gamma ^* \otimes \Lambda ^*$; the associated
product and coproduct are graded commutative in the usual sense. 
\end{exam}

Exponentiality facilitates  calculations of $\hf$:  for $E$ an
exponential functor and $F, G$ functors in $\f$, where $E, F, G$ take
finite-dimensional values, there is a natural isomorphism of vector spaces
\[
\hf (E, F \otimes G) \cong \hf (E, F) \otimes \hf (E, G).
\]
In the graded case, the right hand side has to be treated as a graded 
tensor product. 

\subsection{Endomorphisms of exponential functors}
A graded exponential functor is of finite type if each component is a
finite functor. The endomorphism ring $\fend (E)$ of a graded
exponential functor is a bigraded vector space, which  has additional
structure.

\begin{lem}
Let $E$ be a graded exponential functor of finite type and let ${\fend (E)}
^*$ denote the bigraded dual. Then the following statements hold.
\begin{enumerate}
\item
${\fend (E)}^*$ has the structure of a bialgebra. 
\item
There is a natural coaction 
$\psi \co E \rightarrow E \hat{\otimes} {\fend (E)} ^*$ which satisfies
the following properties: 
\begin{enumerate}
\item
$\psi$ is multiplicative. 
\item
the diagonal $\Delta \co E \rightarrow E \otimes E$ is a morphism of
$\fend (E) ^*$--comodules, where $E\otimes E$ is given the tensor
product comodule structure.
\end{enumerate}
\end{enumerate}
\end{lem}

This result extends to the graded commutative setting.

\subsection{Fundamental examples}

The graded exponential functor $S ^\bullet$ and the bigraded
exponential functor $S ^\bullet \otimes \Lambda ^\bullet$ 
respectively\footnote{The notation
  $E^\bullet$ is used here to avoid confusion with vector space
  duality.} provide the examples of importance to the theory of
unstable modules. 

\begin{exam}
\ 
\begin{enumerate}
\item
The graded exponential functor $S^\bullet$ is of finite type and there is 
a  natural coaction
$$S^\bullet \rightarrow S^\bullet \hctensor {\fend (S^\bullet)}^*.$$
\item
Let $\field$ be a prime field of odd characteristic. The bigraded
exponential functor $\Lambda ^\bullet \otimes \Gamma ^\bullet$ is of 
finite type and there is a natural coaction 
$$S^\bullet \otimes \Lambda ^\bullet
\rightarrow 
(S^\bullet \otimes \Lambda ^\bullet)
\hctensor 
{\fend  (S^\bullet \otimes \Lambda ^\bullet)}^*.$$
\end{enumerate}
\end{exam}

\begin{thm} 
\label{thm:Bend}
For $\field$ a prime field of odd characteristic, there is a natural  
morphism
$\calb \rightarrow \fend(S^\bullet \otimes \Lambda ^\bullet) ^*$ of
$\zed/2$--bialgebras,  which is an isomorphism.
\end{thm}

\begin{proof}
The morphism  exists by the universal property of $\calb$ which is
implicit in the definition. It is straightforward to verify that it is
an isomorphism. 
\end{proof}

\begin{rem}
For $\field$ a prime field of odd characteristic, $\hf
(\Lambda ^m, \Gamma ^n)$ is trivial unless either $m=n=0$ or
$m=n=1$. The vector space $\hf ( \Lambda ^1, \Gamma ^1)$ has dimension
one and the generator $w$ of $\calb$  is dual to a generator of
this vector space.
\end{rem}

\section{Representation categories}
\label{sect-rep-cat}

This section recalls and extends the results of Kuhn \cite{k} on
representation categories which are relevant to the study of the
category of unstable modules. The main result, \fullref{thm:comod_rep}, 
identifies the category $\unst (\calb)$ as a representation category.  The 
functors $r'$, $r_\calb$ introduced in this section are used in the 
analysis of the projective and injective objects of $\unstbigr$ in  
\fullref{section-proj-inj}.

\subsection{Generalities}
Throughout this section, the following hypothesis is supposed to hold
on the pair of categories $\calc, \cals$.

\begin{hyp}
\label{hyp:CS}
The category $\calc$ is abelian and contains  all small inductive
limits, which are exact. The category $\cals \hookrightarrow \calc$ is
a full  small subcategory of $\calc$ with objects $\{ S_i
\}$ indexed over a set $\mathcal{I}$. 
\end{hyp}

Kuhn \cite{k} defines the representation category defined by the category 
$\mathcal{S}$, $\mathrm{Rep} (\mathcal{S} \op)$, as the 
multi-object version of the category of representations (left modules) 
of the ring  $\mathrm{End}(S)\op$, for $S$ an object of $C$. The  
standard example of an object of $\mathrm{Rep} (\mathcal{S} \op)$ is 
given by the $\mathcal{I}$--indexed object 
$ r_{\mathcal{S}} (X)_i  :=  \hom _{\calc} ( S_ i , X)$ for $X$ an object 
of $\calc$. 

There is an adjunction of categories:
\[
l_{\mathrm{S}} \co \mathrm{Rep} (\mathcal{S} \op ) \rightleftarrows 
\mathcal{C} ~{:} r_{\mathcal{S}},
\]
in which the functor $r_{\mathcal{S}}$ is defined as above.

The one-sided Morita equivalence result of Kuhn is basic to the theory:

\begin{thm}[Kuhn {{\cite[Theorem 2.1]{k}}}]
\label{thm-Morita}
The following statements are equivalent.
\begin{enumerate}
\item
$\cals$ generates $\calc$.
\item
$l_{\cals}$ is exact and 
$r_{\cals}$ is fully faithful.
\item
$\calc$ has enough injectives and, for all injectives $I, J$ in $\calc$, 
$r_{\cals}(I), r_{\cals}(J)$ are injective and the functor $r_{\cals}$ 
induces an isomorphism $$
\hom _\calc (I, J) \cong \hom _{\rep{} (\cals \op) } (r_{\cals}(I), 
r_{\cals}(J) ).
$$
\end{enumerate}
Moreover, if these conditions are satisfied, the adjunction counit  
$l_{\cals} r_{\cals} \rightarrow 1_\calc$ is a natural equivalence.
\end{thm}

A set of  projective generators of the category $\rep{}(\cals \op)$ is
given by  $\{ r_{\cals} (S) | S \in \mathrm{Object} (\cals) \}$, 
 by Yoneda's lemma. In the case that $\calc$ is a $k$--linear category 
over a  field $k$, and under a locally finite-type hypothesis, there is a 
dual description of a set of injective cogenerators.

\begin{nota}
\label{nota-rho}
For $\cals, \calc$ as above, where $\calc$ is a $k$--linear category
over  a field $k$,  let 
$\rho \co \calc \rightarrow \rep{} (\cals \op)$ denote the functor 
\[
 X \mapsto \hom _\calc (X, - ) ^*
\]
where $*$ denotes vector space duality and the right hand side is regarded 
as a contravariant functor on $\cals$.
\end{nota}

\begin{prop}
For $\calc$ a $k$--linear category over a field $k$ such that the
vector space $\hom (S, T)$ is of finite dimension, for each pair of
objects $(S, T) $ of $\cals$,  the category $\rep{} (\cals \op)$ has set 
of injective cogenerators $\{ \rho (S)  | S \in \mathrm{Object} ( \cals) 
\}$.
\end{prop}

\begin{proof}
 The proof is straightforward, using vector space duality to  reduce to 
the Yoneda lemma.
\end{proof}

\subsection{Restriction and extension for representation categories}

There are  restriction and extension functors for representation
categories, which are associated to two small subcategories of an
abelian category. This is relevant to the study of the category of
unstable modules over the Steenrod algebra when passing from the
category of objects concentrated in even degree to the full category
of unstable modules (see \fullref{cor:Kan_restrict} below).

Let $\cals, \calc$ satisfy \fullref{hyp:CS} and  suppose moreover 
that there are inclusions of  full subcategories $\cals \subset \calt 
\subset \calc$, where  $\mathcal{T}$ has a set of objects indexed over a 
set $\mathcal{J}$.

There are  abelian representation categories $\rep{} (\cals \op) , \rep{} 
(\calt \op)$ and  canonical adjunctions 
$$l_{\cals} \co \rep{} (\cals \op ) \rightleftarrows \calc ~{:}r_{\cals}
  \qquad\qquad
  l_\calt \co \rep{} (\calt \op ) \rightleftarrows  \calc ~{:}r_\calt,$$
where the functor $r_{\cals}$ is induced by $\hom _\calc ( \cals , -)$ and 
similarly for $r_\calt$.

\begin{prop}
\label{prop:basechange1}
 The inclusion functor $\mathcal{S} \hookrightarrow \mathcal{T}$ induces 
an exact restriction functor  
$
\mathsf{Res} \co
\mathrm{Rep} (\mathcal{T}\op) 
\rightarrow 
\mathrm{Rep} (\mathcal{S} \op). 
$
 Moreover, there is an adjunction 
\[
\mathsf{K} \co \rep{}  (\cals \op) 
\rightleftarrows 
\rep{} (\calt \op) ~{:} \mathsf{Res}
\]
where the functor $\mathsf{Res}$ is defined by restriction and the functor 
$\mathsf{K}$ is induced by Kan extension. The functor $\mathsf{Res}$ is 
exact and the functor $\mathsf{K}$ is right exact.

The Kan functor induces a commutative diagram 
$$
\xymatrix{
\rep{} (\cals \op ) \ar[r] ^{l_{\cals}} 
\ar[d] _{\mathsf{K}}
&
\calc 
\ar@{=}[d]
\\
\rep{} (\calt \op) \ar[r]_{l_\calt } 
&
\calc.
}  
$$
\end{prop}

\begin{proof}
The functor $\mathsf{Res}$ is the evident restriction functor; the functor 
$\mathsf{K}$ is defined by the Kan extension of the functor which 
associates to $r_{\cals} ( S) $ the object $r_\calt (S ) $. 

The commutativity of the given diagram follows from the fact that both 
functors  $\rep{}(\cals \op) \rightarrow \calc$ are left adjoint to the 
functor $r_{\cals} \co \calc \rightarrow \rep{}(\cals \op) $.
\end{proof}

\begin{lem}
Under the hypotheses of \fullref{prop:basechange1}, there is a natural 
transformation
$
l_{\cals} \mathsf{Res} \rightarrow l_\calt
$
of functors from $\rep{}  (\calt \op)$ to $\calc$.
\end{lem}

\begin{proof}
The adjunction unit $1{\rightarrow}r_\calt l_\calt$ induces a natural 
transformation $\mathsf{Res} \rightarrow \mathsf{Res} r_\calt l_\calt$, by 
composition with $\mathsf{Res}$, and the functor $\mathsf{Res} r_\calt 
l_\calt$ is naturally equivalent to $r_{\cals} l_\calt$. The required 
natural transformation is given by adjunction.
\end{proof}

\begin{prop}
\label{prop-base-change}
Under the hypotheses of  \fullref{prop:basechange1},  suppose that the 
objects of $\cals $ generate $\calc$ then the following properties hold.
\begin{enumerate}
\item
The functors $l_{\cals} \co \rep{} (\cals \op) \rightarrow \calc $, 
$l_\calt \co \rep{} (\calt  \op) \rightarrow \calc $ are exact. 
\item 
The natural transformation $l_{\cals} \mathsf{Res } \rightarrow l_\calt$ 
is a natural equivalence. 
\item
The functors $\mathsf{K}, \mathsf{Res}$ induce an equivalence of 
categories 
\[
\rep{} (\cals \op)/ \mathrm{ker}(l_{\cals})  \cong \rep{} ( \calt \op) / 
\mathrm{ker}(l_\calt) .
\]
\end{enumerate}
\end{prop} 

\begin{proof}
The first statement follows from \fullref{thm-Morita}, since the  
hypothesis on $\cals$ implies that the objects of $\calt$ generate $\calc$.

For the second statement, consider the natural transformation $l_{\cals} 
\mathsf{Res} \rightarrow l_\calt$. The functors $l_{\cals} \mathsf{Res}$ 
and $l_\calt$ are exact and send coproducts to coproducts, under the 
hypotheses of the Proposition. Hence, by forming projective resolutions, 
it is sufficient to show that the natural transformation is an equivalence 
 on a set of projective generators of $\rep{} (\calt \op)$. The  Yoneda 
lemma  implies that the objects $r_\calt (T)$, for $T$ objects of $\calt$, 
form such a set of projective  generators. Hence, it is sufficient to show 
that the natural transformation
$ 
l_{\cals} \mathsf{Res} r_\calt \rightarrow l_\calt r_\calt
$ 
induced by composition with $r_\calt$, is an equivalence. Under the 
hypotheses of the Proposition, both functors above are naturally 
equivalent to the identity of $\calc$, by \fullref{thm-Morita}; the 
verification that the above morphism is a natural equivalence is an 
adjunction argument, which is left to the reader.

 The final statement is a corollary of the identification, provided by 
Kuhn, of $\calc $ with the respective localized representation categories. 
\end{proof}

\begin{exam}
The hypotheses of the Proposition are necessary: for example, 
let $\calc$ be the category $\f$ and let $\cals$, $\calt$ be the full
subcategories with sets of objects $\{\Gamma ^1  \}$ and $\{ \Gamma ^n
|n \geq 0\}$  respectively. The category $\rep{} (\cals \op)$ is 
equivalent to the category of vector spaces and it is straightforward to 
see that the natural transformation $l_{\cals} \mathsf{Res} \rightarrow 
l_\calt$ is not a natural equivalence.
\end{exam}

\subsection{Representation categories related to unstable modules}

Kuhn  showed that the category of unstable modules
over the $\field_2$--Steenrod algebra is equivalent to the
representation category associated to the set of objects $\{\Gamma ^* \}$ 
in the category of functors $\f$ defined with respect to $\field_2$--vector
spaces. This result extends to give a description of $\unst'$ in the
case of odd characteristic, which is equivalent to the representation 
category  $\rep{} (\cals \op)$ for the full subcategory $\cals$ of $\f$ 
with 
set of objects $\{\Gamma ^{*} \}$.

This does not extend to a  description of the
category $\unst$; however, the following holds:

\begin{thm}
\label{thm:comod_rep}
$\unst (\calb)$ is equivalent to the representation category 
$\rep{} (\calt \op)$ for the full subcategory $\calt$ of $\f$ with 
set of objects $\{\Gamma ^{*_1} \otimes \Lambda ^{*_2} \}$.
\end{thm}
 
\begin{proof}
This is an immediate consequence of  \fullref{thm:Bend}.
\end{proof}

\begin{cor}
\label{cor-H-injective}
For $n$  a non-negative integer,  the object $ r_\calb (  I_{\field ^n})$ 
of $\unst (\calb)$ is injective. 
\end{cor}

The theorem provides adjunctions 
$$l_\calb\co \unst (\calb )  \rightleftarrows \f ~{:} r_\calb
\qquad\qquad
l' \co \unsteven \rightleftarrows \f ~{:} r '.$$
In terms of the general framework introduced in the previous sections, the 
category $\calc$ corresponds to the  full subcategory $\f_ \omega$ of 
analytic functors in $\f$ (see Henn--Lannes--Schwartz \cite{hls} and
Kuhn \cite{k}).

\begin{nota} For a  full subcategory $\cals \hookrightarrow \calc$ with 
set of objects $\mathbb{S}$, write  $\rep(\{\mathbb{S}\}\op)$ for  the 
associated representation category $\rep (\cals \op) $.
\end{nota} 

\begin{cor}
\label{cor:Kan_restrict}
There is a  diagram of functors, which is  commutative up to natural 
equivalence:
$$
\xymatrix{
\unst (\calb ) \cong 
\rep (\{\Gamma ^* \otimes \Lambda ^* \}\op) 
\ar[rr]^(.7){l_\calb}
\ar@<1ex>[d]^{\mathsf{Res}}
&&
\f
\ar@{=}[d]
\\
\unst ' \cong \rep( \{\Gamma ^* \} \op) 
\ar@<1ex>[u]^{\mathsf{K}}
\ar[rr]_(.7){l'}
&&
\f,
}
$$
in which $l_\calb,l'$ are the exact left adjoints of the adjunction 
between $\f$ and the respective representation categories. Moreover, the 
adjoint functors $(\mathsf{K}, \mathsf{Res}) $ induce an equivalence of 
categories 
$$
\rep (\{\Gamma ^* \otimes \Lambda ^* \}\op)/ \mathrm{ker}(l_\calb)
\cong
\rep (\{\Gamma ^* \}\op)/ \mathrm{ker}(l')
\cong 
\f_\omega
.
$$
\end{cor}

\begin{proof}
The result follows from \fullref{prop-base-change}.
\end{proof}

\part{Further properties}

\section{Projective and injective objects}
\label{section-proj-inj}

The definition of $\unstbigr$ as a category of comodules allows the
 construction of injective cogenerators  and projective
generators in terms of cotensor products. For non-negative integers 
$a,b$, the category $\unstbigr$
contains projective objects $F(a,b)$, analogues of the Massey--Peterson
free unstable modules $F(n) \in \unst$,  and injective objects $J(a,
b)$, analogues of the Brown--Gitler modules $J(n) \in \unst$.  The 
structure of the objects $F(a, b)$ in the category $\unstbigr$  sheds 
light on the structure of the objects $F(n)$ in the category of unstable 
modules. (The reader is referred to Schwartz \cite{s}
for the traditional approach to the projective and injective objects in 
$\unst$).

Throughout this section, let $\field$ be the prime field of
characteristic $p$, where $p$ is an odd prime.  The results  of 
\fullref{section-comodules}, together with the material on graded comodules 
introduced in \fullref{section-comodules-unstable},  will be used in this 
section  without further comment.

%
\subsection[The comodules F(a,b) and J(a,b)]
{The comodules $F(a,b)$ and $J (a,b)$}

Recall that $\vsgrad$ denotes the category of $\zed$--graded vector
spaces and  $\vsbigrad$ denotes the category of
$\zed\times\zed$--graded vector spaces.

Let $\field (n)\in \vsgrad$ denote the graded vector space of total 
dimension one
concentrated in degree $n$ and let $\field (a, b ) \in \vsbigrad $
denote the bigraded vector space of total dimension one concentrated
in bidegree $(a,b)$. There is an exact functor $\vsbigrad \rightarrow
\vsgrad$  which commutes with colimits and which  is therefore
determined by $\field (a, b) \mapsto \field (a+2b)$, for all pairs of 
integers
$(a, b)$.

There is a commutative diagram of exact functors:
$$
\xymatrix{
\unst (\calb ) 
\ar[r]^\Psi
\ar[rd]_{\gamma(\calb)}
&
\unstbigr 
\ar[r]^{\Theta}
\ar[d]^{\gamma (\bbar)}
&
\unst 
\ar[d]^\gamma
&
\unst ''
\ar[l]
\ar[d] ^{\gamma '}
\\
&
\vsbigrad 
\ar[r]
&
\vsgrad
& 
\vsgrad,
\ar@{=}[l]
}
$$
in which the arrows labelled by variants of $\gamma$ are the gradings 
induced by corestriction functors.

\begin{prop}
The functors $\gamma \co \unst \rightarrow \vsgrad$, $\gamma' \co \unst' 
\rightarrow \vsgrad$, $\gamma (\bbar ) \co
\unstbigr \rightarrow \vsbigrad $,  $\gamma (\calb) \co \unst
(\calb) \rightarrow \vsbigrad$ admit both left and right
adjoints. 
\end{prop}

\begin{proof}
The results (cf \fullref{sect:cotensor}) on the existence of right 
and left adjoints to the
corestriction functor extend to these categories $\comodgr$ of graded
right comodules. 
\end{proof}

The following Corollary is immediate.

\begin{cor}
The categories $\unst$, $\unst '$,  $\unstbigr $, $\unst (\calb)$ have 
enough
projective objects and enough injective objects.
\end{cor}

The Brown--Gitler modules and the Massey--Peterson modules have the
following definition from the comodule viewpoint. (cf the material
of  \fullref{sect:cotensor} on cotensor products and duality).

\begin{defn}
For $n$ a non-negative integer, define the following objects of $\unst$, 
considered as the category of right $\atilde$--comodules: 
\begin{enumerate}
\item
$J (n) := \field (n) \Box
_{\field [u]} \atilde$;
\item
$F(n)$ the  right $\atilde$--comodule associated by duality to the left
$\atilde$--comodule $\atilde \Box_{\field [u]} \field (n)$. 
\end{enumerate}
\end{defn}

There is the analogous definition in the bigraded situation: 

\begin{defn}
For  $a, b $ non-negative integers,  define the following
objects of $\unstbigr$, the category of right $\bbar$--comodules:
\begin{enumerate}
\item
$J(a, b) : = \field (a, b) \Box _ {\field [u , \xi_0]} \bbar$; 
\item
$F(a, b)$, the right comodule associated to the left 
$\bbar$--comodule $\bbar \Box_{\field [u , \xi_0]} \field (a,b)$ by
duality.
\end{enumerate}
\end{defn}

The definitions of the objects $F (a, b)$, $J(a,b)$ in terms of adjoints 
to the corestriction functor implies the following characterization, which 
is analogous to the characterization of the unstable modules 
$F(n), J (n)$. 

\begin{prop}
For $a,b$  non-negative integers,   the object $F(a,b)$ is projective in 
$\unstbigr$,  $J (a, b)$ is
 injective in $\unstbigr$ and, for   $M$ an object of $\unstbigr$, there 
are isomorphisms
\begin{enumerate}
\item
$\hom_{\unstbigr} (F (a,b), M) \cong M _{(a,b)}$;
\item
$\hom_{\unstbigr}( M, J (a,b))\cong M _{(a,b)}^*$,
\end{enumerate}
where $M_{(a,b)}$ denotes the homogeneous component of $M$ in bidegree 
$(a,b) $.\end{prop}

It is straightforward to verify the following connectivity result:

\begin{lem}
\label{lem-connectivity}
For $a,b$ non-negative integers,
\begin{enumerate}
\item
$F (a,b) _{(s,t)}= \left\{
\begin{array}{ll}
0 & s<a \text{ or } t<b\\
\field & (s,t)=(a,b).
\end{array}
\right. 
$
\item 
$J (a,b) _{(s,t)}=
\left\{
\begin{array}{ll}
0 & s>a \text{ or } t>b\\
\field & (s,t)=(a,b).
\end{array}
\right. 
$
\end{enumerate}
\end{lem}

The lemma implies that the objects $F (a,b)$ and $J (a,b)$ both have 
fundamental classes in bidegree $(a,b)$, which are unique up to non-zero 
scalar multiple.

\begin{lem}
\label{lem-Theta-canon}
For $a,b$  non-negative integers,
\begin{enumerate}
\item
the object $\Theta J(a,b)$  of $\unst$ is a sub-module of
the Brown--Gitler module $J(a{+}2b)$;
\item
the object $\Theta F (a,b ) $  of $\unst$ is a quotient of the free 
unstable module $F (a{+}2b )$.
\end{enumerate}
\end{lem}

\begin{proof}
Straightforward.
\end{proof}

This implies the following:

\begin{lem}
For $a,b$ non-negative integers, there are isomorphisms:
\begin{enumerate}
\item
$\mathrm{End} _\unst (\Theta J (a, b ) ) \cong \field$;
\item
$\mathrm{End} _\unst (\Theta F (a, b ) ) \cong \field$.
\end{enumerate}
\end{lem}

%
\subsection[Properties of the projective objects of R]
{Properties of the projective objects of $\unstbigr$}

The projective $F(a,b)$ is
defined as the left comodule $\bbar \Box_{\field [u, \xi_0]}
\field (a,b)$ and hence its underlying vector space depends only upon the 
right $\field [u, \xi_0]$--comodule structure of $\bbar$. The following 
result reduces the study of these projective generators, via tensor 
products, to the cases where either $a$ or $b$ is zero.

\begin{prop}
\label{prop-bigrad-tensor}
For $a,b$ non-negative integers, there is an
isomorphism of objects $F (a, b) \cong F (a, 0) \otimes F (  0, b )$ in 
$\unstbigr$. 
\end{prop}

Recall that $\field [\xi_j |j \geq 0]$ has a natural  bialgebra
structure and that there exists  a morphism of bialgebras $\field [\xi_j 
|j \geq 0] \rightarrow \field [\xi_0]$,  which induces a right $\field 
[\xi_0]$--comodule structure upon $\field [\xi_j |j \geq 0]$. The 
$\zed/2$--graded algebra $\field [u ] \otimes \Lambda ( \tau_i |i \geq 0) 
$ has the structure of a right $\field [u] $--comodule with respect to the 
multiplicative  structure morphism which is induced by
 $
u \mapsto  u \otimes u, 
\tau_i \mapsto  \tau_i \otimes u .
$

\fullref{prop-bigrad-tensor} is a formal consequence of the
following observation, using the properties of cotensor products over
a tensor product of coalgebras.

\begin{lem}
As a right $\field[u, \xi_0]$--comodule, $\bbar$ is isomorphic to the 
exterior tensor product of the right $\field [\xi_0]$--comodule $\field 
[\xi_j |j \geq 0]$ and of the right $\field [u]$--comodule $\field [u] 
\otimes \Lambda (\tau_i)$. 
\end{lem}

\begin{proof}
Immediate.
\end{proof}

Recall that there is a commutative diagram of functors
$$
\xymatrix{
\f 
\ar[rr]^(.3){r_\calb} 
\ar[rrd]_{\Psi r_\calb}
&&
\rep  (\{\Gamma^* \otimes \Lambda ^*\} \op) \cong \unst (\calb)
\ar[d] ^{\Psi}
\\
&&\unstbigr
\\
\f
\ar[rru]^{\hat{\calo}r'}
 \ar[rr]_(.3){r'}
&&
\rep (\{\Gamma ^* \} \op)\cong \unsteven
\ar[u]_{\hat{\calo}},
}
$$
in which $\Psi$  and $\hat {\calo}$ denote the exact corestriction 
functors.  

\begin{prop}
\label{prop:identify-via-functors}
For $a, b$ non-negative integers, there are natural isomorphisms in 
$\unstbigr$
\begin{enumerate}
\item
$F (0, b) \cong \hat {\calo} r' (\Gamma ^b)$ 
\item
$F (a, 0) \cong \Psi r_\calb (\Lambda ^a) $.
\end{enumerate}
\end{prop}

\begin{proof}
(Indications) \qua
 It is straightforward to show that there is a surjection  $F (0, b) 
\rightarrow \hat{\calo} r ' (\Gamma ^b)$, using the defining property of 
$F (0, b)$. This can be seen to be an isomorphism by comparing Poincar{\'e}
series. 
The second statement admits a similar proof; namely, there is a surjection 
$F (a, 0 ) \rightarrow \Psi r _\calb (\Lambda ^a)$ and the result follows 
by comparing Poincar{\'e} series.    
\end{proof}

\begin{rem}
\ 
\begin{enumerate}
\item
It is not true in general that there is a surjection $F (0, b ) 
\rightarrow \Psi r (\Gamma ^b)$; the obstruction is the isomorphism 
$\Lambda ^1 \rightarrow \Gamma ^1$. 
\item
The structure of the objects $F (0, b)$ is analogous to the structure
of the projective generators $F (b)$ in the category of unstable
modules at the prime two (see \cite[Proposition 1.6.3 and Proposition 
1.7.3]{s}).
\end{enumerate}
\end{rem}

\begin{lem}
For  a positive integer $a$,  there is a non-trivial surjection $F (a, 0 ) 
\rightarrow \Sigma^{(1,0)} F (a- 1, 0)$, which is unique up to non-trivial 
scalar multiple.
\end{lem}

\begin{proof}
Straightforward.
\end{proof}

The description of the objects $F (a,0)$ given in 
\fullref{prop:identify-via-functors} leads to the direct proof of
the following Proposition.

\begin{prop}
\label{prop-g-filtration}
For  a positive integer $a$, there exists a finite filtration 
$$
0 = g_{-1} F( a, 0) \subset g_0 F (a, 0 ) \subset \ldots g_a F (a, 0 )
= F (a, 0 ), 
$$
such that the filtration quotients are identified, for $0 \leq j \leq a$, 
by:
\[
g_{j} F(a,0) / g_{j-1} F (a,0) 
\cong 
\Sigma ^{(j, 0 )} \hat{\calo} r ' (\Lambda ^{a-j}) .
\]
In particular, there is a monomorphism 
$
\hat{\calo} r ' (\Lambda ^a ) \hookrightarrow F (a,0),
$
 which fits into a short exact sequence in $\unstbigr$:
\[
0\rightarrow 
\hat{\calo} r ' (\Lambda ^a ) \rightarrow F (a,0)
\rightarrow 
\Sigma^{(1, 0)} F (a-1, 0)
\rightarrow 0 .
\]
\end{prop} 

\begin{proof}
(Indications) \qua
The filtration can be deduced by using the exponential property of the
graded functor $\Lambda ^*$ in $\f$ and the fact that $\Lambda ^n$ is
a simple object for each integer $n$. For the final statement, it is 
straightforward to identify the kernel of the morphism $F (a, 0) 
\rightarrow \Sigma ^{(1,0)} F (a-1, 0)$. 
\end{proof}

\begin{nota}
For a non-negative integer $a$, define  $\Phi F (a, 0)$ to  be $F
(a,0) $ for $a \in \{0,  1\}$ and, for $a >1$,  via the short exact 
sequence in $\unstbigr$:
\[
0
\rightarrow 
\Phi F (a, 0) 
\rightarrow 
F (a, 0) 
\rightarrow 
\Sigma^{(2, 0)} F (a- 2, 0)
\rightarrow 0
. 
\]
(The convention is adopted that $\Phi$ is associated to the double
suspension $\Sigma^{(2,0)}$ and to the single suspension
$\Sigma^{(0,1)})$.
\end{nota}

\begin{lem}
\label{lem-decomp-Phi}
For $a >1$ an integer, there is a  short exact sequence in $\unstbigr$:
\[
0 \rightarrow
\hat{\calo} r ' ( {\Lambda}^a)
\rightarrow 
\Phi F (a, 0) 
\rightarrow 
\Sigma^{ ( 1,0)}
\hat{\calo} r '
( {\Lambda}^{a-1})
\rightarrow 0.
\]
\end{lem}

\begin{proof}
The result follows  from the final statement of 
\fullref{prop-g-filtration}.
\end{proof}

\subsection
[Properties of the projective objects in t]
{Properties of the projective objects in $\unst$}

The results of the previous section allow the description of the
projective generators of $\unst$ in terms of the forgetful functor
$\Theta \co \unstbigr \rightarrow \unst$ and the projective objects
$F (a, 0) , F (0, b )  \in \unstbigr$. \fullref{prop-fn-equalizer} 
expresses $F(n)$ in terms of an equalizer diagram; the result is 
essentially a formal consequence of the definition of the corestriction 
functor $\unstbigr \rightarrow \unst$. 

This analysis is used to provide a filtration of the objects $F(n) $, 
which is given in \fullref{thm-filt-quot-fn}; the associated graded has an 
explicit description which can be regarded as being a natural extension of 
the analysis of the structure of the objects $F(n) $ at the prime two.

\begin{lem}
For $n$ a non-negative integer and  $a, b $ non-negative integers such 
that $n = a + 2b $, there is a  
non-trivial morphism 
 $
\mu_{a, b} \co 
F (n ) \rightarrow \Theta F (a, b), 
$  
unique up to non-zero scalar multiple, which is surjective. 
\end{lem}

\begin{proof}
The existence of the surjection is given by \fullref{lem-Theta-canon}; the 
fact that it is unique up to non-zero scalar multiple follows from the 
fact that $\Theta F (a,b)$ is of dimension one in degree $a+2b$, which 
implies that $\mathrm{Hom}_\unst (F (n) , \Theta F( a, b) ) \cong  \field$.
\end{proof}

Recall that there are surjections in $\unstbigr$,  for integers $a \geq 2$ 
and $b \geq 0$:
\begin{eqnarray*}
F (a, 0 ) &\rightarrow &\Sigma^{(2,0)}  F (a-2, 0)
\\
F (0, b+1) &\rightarrow & \Sigma ^{(0, 1)} F ( 0, b).
\end{eqnarray*}
These induce surjections 
\begin{align*}
l_{a, b} \co &\Theta F (a, b)  \rightarrow  \Sigma ^2 \Theta F (a-2, b) \\
r_{a-2, b+1} \co &\Theta F (a-2, b+1) \rightarrow  \Sigma ^2 \Theta F (a-2, 
b)
\end{align*}
in the category $\unst$.

\begin{lem}
For integers $a \geq 2, b \geq 1$ satisfying  $a + 2b= n$,  the following 
diagram commutes up to non-zero scalar multiple
$$
\xymatrix{
F (n) 
\ar[d]_{\mu_{a, b }}
\ar[rr]^(.4){\mu _{a-2, b+1}}
&&
\Theta F (a-2, b+1) 
\ar[d]^{r_{a-2, b+1}}
\\
 \Theta F (a, b) 
\ar[rr]_(.4){l_{a, b }}
&&
\Sigma ^2 \Theta F (a-2,b).
}
$$ 
\end{lem}

\begin{proof}
This result follows immediately from  \fullref{lem-connectivity}.
\end{proof}

\begin{prop}
\label{prop-fn-equalizer}
For  $n$  a non-negative integer,  the object $F (n)$ identifies with the 
equalizer of the diagram 
$$
\xymatrix{
\bigoplus _{a + 2b =n} \Theta F (a, b) 
\ar@<.5ex>[r]^{l_{a,b}} 
\ar@<-.5ex>[r]_{r_{a,b}} 
&
\bigoplus _{c + 2d =n} \Theta F (c, d)
.
}
$$
In particular, the morphism
$$
\bigoplus_{a+2b=n} \mu_{a, b} \co 
F(n) 
\rightarrow 
\bigoplus_{a+2b=n}
\Theta F (a, b )
$$
is a monomorphism.
\end{prop}

The Proposition is a formal consequence of the definition of the functor 
$\unstbigr \rightarrow \unst$; it can be proved explicitly using the 
identification provided  by \fullref{lem-surjectivity-result} below. 
Recall that the bialgebra $\bbar$ has the structure of a right
$\field[u]$--comodule via the corestriction associated to the morphism
of Hopf algebras $\field [u,\xi_0] \rightarrow \field[u]$ given by $\xi_0 
\mapsto u ^2$.

\begin{lem}
\label{lem-surjectivity-result}
For  $n$  a non-negative integer,
\begin{enumerate}
\item
 the surjection
$\bbar \rightarrow \atilde$ induces a surjection 
$
\bbar \Box_{\field [u]} \field (n) 
\twoheadrightarrow
\atilde \Box_{\field [u]} \field (n), 
 $
 which is a morphism of left $\atilde$--comodules, where $\bbar
\Box_{\field [u]} \field (n) $ is given the corestricted structure;
\item
the left $\bbar$--comodule $\bbar \Box_{\field [u]} \field
(n) $ is isomorphic to $\bigoplus_{a+2b = n} \bbar  \Box_{\field [u,
\xi_0]} \field (a, b)$.
\end{enumerate}
\end{lem}

\begin{proof}
The surjectivity of the morphism in (1) follows from the right exactness
of $- \Box_{\field [u]} \field (n)$, which reflects a coflatness
property. The proof of the isomorphism of (2) is straightforward.
\end{proof}

\fullref{prop-fn-equalizer} gives rise to a  filtration of the
object $F(n)$ in the category $\unst$ as follows. The choice of the
filtration is motivated by the consideration of the structure of the
object $F (2)$ in the category $\unst$, for which there is a short exact 
sequence:
\[
0 \rightarrow
\Theta \Phi F (2, 0) 
\rightarrow 
F (2) 
\rightarrow 
\Theta F (0, 1) 
\rightarrow 
0
.
\]

\begin{nota}
For $j \geq -1$ an integer, let $f_j F (n) $ denote the subobject of $F(n) 
 \in \unst$ defined by the  kernel of the morphism
\[
F(n) 
\rightarrow 
\bigoplus 
_{\substack{a+2b = n\\ b \geq j+1} }
\Theta F (a, b)
.
\]
This defines an increasing filtration of $F(n)$, with  $f_{-1} F (n)$  
zero and $f_{j} F(n) = F(n)$ for $2j \geq n$. 
\end{nota}

The filtration quotients can be identified explicitly as follows. In the 
Lemma below, the morphism $l _{a,b} $ is taken to be zero if the integer 
$a$ is in $\{0, 1\}$.

\begin{lem}
\label{lem:l-kernel}
For $a, b \geq 0$ integers, the kernel of the morphism $l_{a,b} \co \Theta F 
(a, b ) \rightarrow \Sigma ^2 \Theta F (a-2, b) $ identifies with the 
object 
 $
\Theta \{\Phi F (a, 0 ) \otimes F (0, b)\}.
$ 
\end{lem}

\begin{proof}
The lemma  follows immediately from the definition of the objects $\Phi F 
(a,0)$. 
\end{proof}

\begin{lem}
\label{lem-mono-filt-quot}
For $a,b$ non-negative integers such that $a+2b =n$, the morphism 
$\mu_{a,b}$ induces a monomorphism
\[
f_b 
F (n) / f_{b-1} F (n) 
\hookrightarrow 
\Theta 
\{\Phi F (a, 0 ) \otimes F (0, b)\}.
\] 
\end{lem}

\begin{proof}
This is an immediate consequence of the definition of the filtration and of 
\fullref{lem:l-kernel}.
\end{proof}

\begin{thm}
\label{thm-filt-quot-fn}
For $n$ a positive integer, the object $F(n)\in \unst$ has a finite
increasing filtration $\{ f_b F(n) \}$   such that the filtration 
quotients are of the form
\[
f_b 
F (n) / f_{b-1} F (n) 
\cong 
\Theta \{
\Phi F (a, 0 )  \otimes F (0,b)\}, 
\] 
where $a,b$ are non-negative integers such that $a+2b = n$. 
\end{thm}

\begin{proof} (Indications) \qua 
The theorem is proved by showing that the monomorphisms defined in 
\fullref{lem-mono-filt-quot} are isomorphisms. This can be proved  by 
analysing the monomial basis of $F(n)$ in conjunction with 
\fullref{prop-fn-equalizer}; indeed, it is sufficient to use a comparison 
of Poincar{\'e} series.  
\end{proof}

A  related analysis of the structure of $F (n)$ is given in Schwartz
\cite{sch}, 
from a  different viewpoint. The advantage of the approach given above is 
that it leads immediately to the following description of the filtration 
quotients in a finite filtration of $F (n)$.

\begin{cor}
For $n$ a positive integer, $F(n)$ has a finite filtration with  
associated graded
\[
\bigoplus _{\substack{a+2b=n\\a \geq 0}}
\calo r ' (\Lambda ^a \otimes \Gamma ^b) 
\oplus
\bigoplus
_{\substack{a+2b=n\\
a \geq 1} }
\Sigma \calo r' ( \Lambda ^{a-1} \otimes \Gamma ^b ).  
\]
  \end{cor}  

\begin{proof}
The Corollary follows from an analysis of the unstable modules $\Theta
(\Phi F(a, 0 ) \otimes F (0, b) )$, by \fullref{thm-filt-quot-fn}. 
\fullref{lem-decomp-Phi} implies that there
is a short exact sequence
\[
0 \rightarrow
\hat{\calo} r ' ( {\Lambda}^a) \otimes F(0, b) 
\rightarrow 
\Phi F (a, 0) \otimes F (0, b)
\rightarrow 
\Sigma^{ ( 1,0)}
\hat{\calo} r '
( {\Lambda}^{a-1}) \otimes F (0, b)
\rightarrow 0
\]
and $F(0, b)$ identifies with $\hat {\calo} r ' (\Gamma ^b) $, by
\fullref{prop:identify-via-functors}. Moreover, exponentiality
implies that there is an identification of  $\hat{\calo} r' (\Lambda ^a)
\otimes \hat{\calo} r' (\Gamma ^b)$ with $\hat{\calo} r' (\Lambda ^a
\otimes \Gamma ^b)$ and likewise for the term involving $a-1$. The
result follows by applying the functor $\Theta$, using the
identification $\calo = \Theta \hat{\calo}$. 
\end{proof}

\subsection
[The injective cogenerators of R and t]
{The injective cogenerators of $\unstbigr$ and $\unst$}
\label{subsect-cogen-inj}

In this section, the injective cogenerators $J (a, b)$ are
analysed. It is shown that the classical results concerning the
structure of the Brown--Gitler modules over an odd prime $p$ arise from
the category $\unstbigr$. In particular, a description of the injective 
cogenerators is given in terms of  the functor 
$\Psi \rho \co \f \rightarrow \unstbigr$. 

\begin{prop}
\label{prop-J-tensor}
For $a,b$ non-negative integers,  there are canonical isomorphisms in 
$\unstbigr$:
\begin{enumerate}
\item
$J (a, 0) \cong \field (a, 0) $
\item
$J (a, b) \cong J (a,0) \otimes J (0, b)$.
\end{enumerate}
\end{prop}

\begin{proof}
The result is  a formal consequence of the observation that, as a left
$\field [u,  \xi_0]$--comodule, the bialgebra $\bbar$ is the  exterior 
tensor product of the algebra $\field [u]$, considered as a left 
$\field[u]$--comodule, and the algebra $\field [\xi_j |j \geq 0]\otimes 
\Lambda (\tau_i |i \geq 0)$, which has the structure of a left 
$\field[\xi_0]$--comodule.
\end{proof}

The Proposition implies that the problem of understanding the
structure of the  objects $J (a, b)$ can be reduced to understanding the 
objects $J (0, b)$ in $\unstbigr$, up to suspension. 

\begin{lem}
For $b$ a non-negative integer, there is a canonical monomorphism 
$
\Sigma^{(0, 1)} J (0, b) \hookrightarrow 
J (0, b+1)
$
 in the category $\unstbigr$, unique up to non-zero scalar multiplication.
\end{lem}

\begin{proof}
Straightforward.
\end{proof}

Recall that the functor $\rho \co \f \rightarrow \rep (\{\Gamma ^* \otimes 
\Lambda ^* \}\op ) \cong \unst (\calb)$ is defined by $F \mapsto \hf (F, 
-) ^*$. There is  the composite functor
\[
\Psi \rho \co \f  \rightarrow \unstbigr,
\]
where $\Psi$ is the forgetful functor $\unst (\calb) \rightarrow 
\unstbigr$. The bigraded  exponential property of the functors $\Gamma ^* 
\otimes \Lambda ^*$ implies the following variant of a  standard result. 

\begin{prop}
For $F, G$ objects of $\f$,  there is a natural isomorphism in 
$\unstbigr$: 
\[
\Psi \rho (F \otimes G ) 
\cong 
\Psi \rho (F) \otimes \Psi \rho (G),
\]
where the tensor product on the right hand side denotes the bigraded 
tensor product.
\end{prop}

The following result gives an elegant description of the objects 
$J (0, n)$.

\begin{prop}
\label{prop-J0n}
For $n$  a non-negative integer, there is an isomorphism 
$
J (0, n) 
\cong 
\Psi \rho (\Gamma ^n)$ 
in $\unstbigr$. In particular, in bidegree $(a, b)$, there is an 
isomorphism of vector spaces 
\[
J (0, n) _{a, b } \cong \hf (\Gamma ^n, \Lambda ^a \otimes \Gamma ^b ) ^*.
\]
\end{prop}

\begin{proof} (Indications) \qua
The result follows from  the identification of the category $\unst 
(\calb)$ as the representation category $\rep (\{\Gamma ^*\otimes \Lambda 
^* \}\op) $, the definition of the functor $\Psi$, together with the 
observation that the morphism $\Lambda ^1 \rightarrow \Gamma ^1$ does not 
intervene in the calculation of $\hf (\Gamma ^n , \Lambda ^* \otimes 
\Gamma ^*)$.
\end{proof}

\begin{cor}
There is a  surjection $J (0, 1) \stackrel{\epsilon(0,1)} {\longrightarrow} 
\field (1, 0) $, which is
unique up to non-zero scalar multiple.
\end{cor}

\begin{proof}
\fullref{prop-J-tensor} implies that $\field (1,0) \cong J
(1,0) $, hence it is sufficient to show that the object $J (0, 1) $
has dimension one in bidegree $(1, 0)$, using the representing
property of $J(1,0)$; this follows from \fullref{prop-J0n}.
\end{proof}

\begin{rem}
The analogous description does not hold for $J (a,0)$, in general; for  
$\hf (\Lambda ^a, \Gamma ^* \otimes \Lambda ^*)$ has total dimension two, 
for a positive integer $a $, since the diagonal morphism $\Lambda 
^{a}\rightarrow \Lambda ^1 \otimes \Lambda ^{a-1}$ induces a morphism 
$\Lambda ^{a}\rightarrow \Gamma ^1 \otimes \Lambda ^{a-1}$, using the 
isomorphism $\Lambda ^1 \cong\Gamma ^1$.
\end{rem}

This result gives the following information on the structure of the 
objects of the form $\Psi \rho (F)$ in the category $\unstbigr$.

\begin{cor}
For  $F \in \f$ a finite functor, there exists an exact sequence in 
$\unstbigr$
\[
\bigoplus _j J (0, c_j ) \rightarrow 
 \bigoplus _{i} J (0, b_i) \twoheadrightarrow \Psi \rho (F) \rightarrow 0
\]
for finite sets of non-negative integers $\{b_i \}$, $\{c_j \}$.
\end{cor}

\begin{proof}
 \fullref{thm-emb-thm}  implies that there exist  finite sets of integers 
$\{b_i \}$, $\{c_j \}$  and an exact sequence  $ \bigoplus _j \Gamma 
^{c_j}  \rightarrow  \bigoplus _i \Gamma ^{b_i} \rightarrow F \rightarrow 
0$. The result follows since $\rho$ is  right exact.  
\end{proof}

\begin{nota}
Let $\rho ' \co \f \rightarrow \rep (\{\Gamma ^* \}\op )$ denote the functor 
$\rho$ associated to $\{\Gamma^* \}$, as in \fullref{nota-rho}.
\end{nota}

Recall that  there is a  natural embedding $\unsteven \hookrightarrow
\unstbigr$ which is denoted  by $\hat{\calo}$.  There is a natural 
embedding 
$
\hat{\calo} \rho' (\Gamma ^n) 
\hookrightarrow
J (0, n)
$ in $\unstbigr$,  which corresponds to the beginning of a filtration
with  associated graded described explicitly by the following result.

\begin{prop}
For a non-negative integer $n$, the object $J (0,n)$ has a finite
filtration with  associated graded 
\[
\bigoplus
\Sigma ^{(d,0)} \hat{\calo} \rho'( \Gamma ^{n-d}), 
\]
where $d$ ranges over the set of integers $0 \leq d \leq n$ which can be 
expressed as a sum of at most $n$ pairwise distinct powers of $p$. 
\end{prop}

\begin{proof}
The result is deduced by using the exponential property of $\Gamma^* $
to obtain a direct sum decomposition at the level of vector spaces, 
which is the associated graded to a filtration in $\unstbigr$.
\end{proof}

There are analogues of the Mahowald exact sequences in the category
$\unstbigr$ (cf \cite[Proposition 2.3.4]{s}). To state the result,
define the following morphisms which are induced by the Verschiebung
$\Gamma^{np}\smash{\stackrel{\mathcal{V}}{\rightarrow}} \Gamma ^n$. (Recall
that the Verschiebung is the morphism which is dual to the Frobenius $p$th 
power morphism).

\begin{defn}
For $n$ a positive integer, 
\begin{enumerate}
\item
let $\mathcal{V}_n \co J (0, np) \rightarrow J (0, n) $ denote the
surjective morphism which is induced under the functor $\Psi \rho \co \f
\rightarrow \unstbigr$ by the Verschiebung 
$\Gamma ^{np} \smash{\stackrel{\mathcal{V}}{\rightarrow}}
\Gamma ^n$; 
\item
let $\mathcal{V}_n' \co J (0,np+1) \rightarrow \field (1,0) \otimes J(0,
n)\cong \Sigma^{(1,0)}J(0,n)$ denote the composite morphism
\[
J(0, np+1) \rightarrow J(0, 1) \otimes J (0,n)
\stackrel{\epsilon(0,1) \otimes 1}{\longrightarrow}  \field (1,0) \otimes 
J(0,n)
\]
where the first morphism is induced under $\Psi \rho$ by  the composite
\[
\Gamma ^{np+1} \stackrel{\Delta}{\rightarrow} \Gamma^1 \otimes \Gamma 
^{np} 
 \stackrel{1 \otimes\mathcal{V}}{\rightarrow}   \Gamma^1 \otimes \Gamma 
^{n},
\]
in which $\Delta$ denotes the coproduct. 
\end{enumerate}
\end{defn}

\begin{prop}
\label{prop-Mahowald-sequences}
For  $n$ a positive integer,  there are short exact sequences in 
$\unstbigr$:
\begin{enumerate}
\item
$
0 \rightarrow \Sigma^{(0, 1)} J (0, np -1) 
\rightarrow 
J (0, np ) 
\stackrel{\mathcal{V}_n} {\rightarrow }
J (0, n) 
\rightarrow 0;
$
\item
$
0 
\rightarrow
\Sigma^{(0, 1)} J (0, np ) 
\rightarrow 
J (0, np +1) 
\stackrel{\mathcal{V}'_n} {\rightarrow }
\Sigma^{(1,0)}
J (0, n)
\rightarrow 0. 
$
\end{enumerate}
If $m$ is a positive  integer such that $m \not \equiv 0, 1 \mod {p}$, 
then the  canonical monomorphism 
\[
\Sigma^{(0, 1)}J (0, m-1 )
\cong 
J (0, m) 
\]
is an isomorphism.
\end{prop}

\begin{proof} (Indications) \qua
The result can be proved by a calculation in terms of the description of 
$J (0, *) $ as $\Psi \rho (\Gamma ^*)$. 
\end{proof}

\begin{rem}
The consideration of the bigraded situation naturally gives rise to
the suspension which appears in the second Mahowald exact sequence
(see  \cite[Proposition 2.3.4]{s}), which  corresponds to the presence
of  the Bockstein operation.
\end{rem}

\begin{prop}
\label{prop-BG}
For $n$  a positive integer, there are canonical isomorphisms
\begin{enumerate}
\item
$J (2n ) \cong \Theta J (0, n)$
\item
$J (2n+1) \cong \Theta J (1, n)$,
\end{enumerate}
where $\Theta \co \unstbigr \rightarrow \unst$ denotes the canonical 
forgetful functor.
\end{prop}

\begin{proof} (Indications) \qua
There are many ways of approaching the proof of this result. For example, 
use the canonical (up to non-zero scalar multiple) monomorphisms
$\Theta\co J(a,b) \hookrightarrow J (a+2b)$. It is straightforward to use the 
Mahowald exact sequences to compare Poincar{\'e} series in the relevant cases 
of the Proposition.
\end{proof}

\section{New proofs of fundamental results}
\label{sect:new-proofs}

The analysis of the categories $\unst (\calb)$, $\unstbigr$ in
relation to the category $\unst$ given in this paper yields 
direct proofs of the foundational results of the theory of
nil-localization of the category of unstable modules in odd
characteristic, generalizing the approach of Kuhn available over the
field $\field_2$. 

\subsection
[Injectivity of H*(Bn) in t]
{Injectivity of $H^* (B \field ^n)$ in $\unst$}

This section is devoted to giving a self-contained  proof of  the
injectivity of the object $H^*(B\field^n)\cong \tilde{\Theta} r_\calb
(I_{\field ^n}) $ of $\unst$ \fullref{thm-Morita}  applied to the 
representation category $\rep(\{\Gamma ^* \otimes \Lambda ^* \} \op)$. The 
proof relies on the analysis of the injective cogenerators of the 
categories $\unst, \unstbigr$  in \fullref{section-proj-inj}.

\begin{prop}
\label{prop-tilde-Theta-inj}
Let $a$ be a non-negative integer, then the object $\tilde{\Theta} \rho 
(\Gamma ^a) $ is injective in $\unst$.
\end{prop}

\begin{proof}
 The injective cogenerators of $\unst$ are considered in 
\fullref{subsect-cogen-inj}; in particular, combining \fullref{prop-J0n} 
and \fullref{prop-BG} implies that the object $\tilde{\Theta} \rho (\Gamma 
^{a_i})$ identifies with an injective object in $\unst$.
\end{proof}

\begin{lem}
Let $n$ be a non-negative integer, then there exists a monomorphism in 
$\unst (\calb)$ 
\[
r_\calb (I_{\field ^n}) \cong H ^{\otimes n} 
\hookrightarrow 
\prod_{i\geq 0} \rho ( \Gamma ^{a_i} ),
\]
where the sequence of integers $a_i$ can be taken to have limit $\infty$ 
as $i$ goes to $\infty$.
\end{lem}

\begin{proof}
The structure theory of $\unsteven \cong \rep(\{\Gamma ^* \}\op) $ implies 
that  there exists a morphism   $\alpha \co r_\calb (I_{\field ^n}) \cong H 
^{\otimes n} 
\rightarrow 
\prod_{i\geq 0} \rho ( \Gamma ^{a_i} )$ in $\unst (\calb)$ such that the 
restriction $\mathsf{Res} (\alpha) $ in $\unst'$ is a monomorphism. 

\fullref{thm-emb-thm} implies that $\alpha$ is injective, as
follows; suppose that $(a, b)$ is a pair of non-negative integers,
then  \fullref{lem-emb-tensor} implies readily  that there is a surjection 
$
\phi \co \Gamma ^T \rightarrow \Lambda^a \otimes \Gamma ^b$, for an
integer $T$.  The morphism $\alpha$ is injective in bidegree $(0, T)$,
by hypothesis; the  surjection $\phi$ implies that the morphism $\alpha$ 
is injective in bidegree $(a,b)$.
\end{proof}

\begin{thm}
\label{thm:Hinj}
For a non-negative integer $n$, the object $H^* (B \field^n)$ is injective 
in $\unst$. 
\end{thm} 

\begin{proof}
There is a monomorphism 
$
r_\calb (I_{\field ^n}) \cong H ^{\otimes n} 
\hookrightarrow 
\prod_{i\geq 0} \rho ( \Gamma ^{a_i} )
$
in $\unst(\calb)$, by the previous Lemma.  The injectivity of $r_\calb 
(I_{\field ^n})$ in $\unst (\calb)$ (\fullref{cor-H-injective}) yields a 
retraction. 

Applying the functor $\tilde{\Theta} \co \unst (\calb) \rightarrow \unst$ 
yields a split monomorphism in $\unst$. The functor $\tilde{\Theta}$ 
commutes with products and  the objects $\tilde{\Theta} \rho (\Gamma 
^{a_i})$ are injective in $\unst$, by \fullref{prop-tilde-Theta-inj}, 
hence the result follows.
\end{proof}

\subsection{Nilclosure for unstable modules}

This section  makes explicit the structure of a nil-closed unstable
module, when the field $p$ has odd characteristic.  \fullref{thm:Hinj} 
implies that there is an exact functor $l \co \unst \rightarrow \f$ and an 
adjunction
\[
l\co \unst \rightleftarrows \f ~{:} r,
\]
 where the functor  $l \co \unst \rightarrow \f$ is defined to be  the 
functor which sends the object $M$ of $\unst$ to the functor  $\field ^n 
\mapsto \hom_\unst (M , H^* (B\field ^n) )  '$, where $'$ denotes the 
profinite dual.

\begin{nota}
Let $\tilde{\Theta} \co \unst (\calb) \rightarrow \unst$ denote the exact 
corestriction functor, which identifies with the composite $\unst (\calb) 
\rightarrow \unstbigr \rightarrow \unst$. 
\end{nota}

The definition of the categories as comodule categories  implies the 
following result. 

\begin{lem}
The functor $\tilde{\Theta} \co \unst(\calb) \rightarrow \unst$ admits a 
left adjoint $\lambda \co \unst \rightarrow \unst (\calb) $. 
\end{lem}

The object $H^* (B \field ^n)$ of $\unst$ identifies with the object
$\tilde{\Theta} r_\calb I_{\field^n}$, hence the following result
follows, using the identification of the functor $l_\calb$ which is
provided by the theory of \fullref{sect-rep-cat}.

\begin{cor}
\label{cor:adjunct-factor}
The adjunction $l \co \unst \rightleftarrows \f ~{:} r $ is the composite of 
the adjunctions $\lambda \co \unst \rightleftarrows \unst (\calb ) ~{:} 
\tilde{\Theta}$ and $l_\calb \co \unst (\calb ) \rightleftarrows \f ~{:}
r_\calb$. 
\end{cor}

\begin{rem}
The theory of \cite{hls} includes the identification of the kernel of
the functor $l$ as the full subcategory of nilpotent unstable  modules, 
which is omitted above.
\end{rem}

Using the terminology of localization of abelian categories as in
\cite{hls}, \fullref{cor:adjunct-factor} implies the following result.

\begin{cor} 
\label{cor:nil-closed}
An object of $\unst$ is nil-closed if and only it is isomorphic to an 
object of the form $\tilde{\Theta} r_\calb (F)$, for $F $ an object of 
$\f$.
\end{cor}

\begin{exam}
The object $\Theta F (a,0)$ is a nil-closed object in $\unst$, whereas 
$\Theta F (0, b)$ is not nil-closed, for $b \geq 1$. It follows that the 
object $\Theta F (a,b) $ is not nil-closed in $\unst$,  for $b \geq 1$. 
There is a canonical  monomorphism
 $
\Theta F (a, b )  \hookrightarrow \tilde{\Theta} r_\calb (\Lambda ^a 
\otimes \Gamma ^b)
 $,
which represents the nil-closure, where $\tilde{\Theta} \co \unst (\calb ) 
\rightarrow \unst $ is the corestriction functor.
\end{exam}

\part[Motivic vistas at the prime 2]{Motivic vistas at the prime $2$}

\section{Bigraded and ordinary unstable modules at the prime two}
\label{sect-motivic}

The construction of a category of bigraded unstable modules
$\unstbigr$ at the prime two proceeds as in the odd characteristic
case by considering the endomorphisms of the primitively-generated
Hopf algebra $\field [y]/ y^2 \otimes \field [x]$, where $\field$ is
taken to denote $\field_2$ throughout this section.

There is no direct analogue of
the exact functor $\Theta$ to the category of unstable modules; the 
categories
are related through a category $\unstmot$,
which is a perturbation of the category of unstable modules;  the
category $\unstmot$  can be  related to modules over the motivic Steenrod
algebra (Voevodsky \cite{v1,v2}). 
In particular, motivic cohomology is a bigraded
cohomology theory, hence the presence of the bigrading is essential.

This section outlines the construction of $\unstmot$ and its relation
to $\unst$; some of the results are implicit in the work of Yagita 
\cite{y1,y2}.

\subsection[The category M]{The category $\unstmot$}
The underlying category used in this section is the category of modules
over the polynomial algebra $\field [\tau]$. In the study of motivic
cohomology, this corresponds to the fact that the coefficient ring is
not concentrated in a single degree.

The category $\unstmot$ is constructed by considering endomorphisms of
the algebra $\hmot$ which is defined as follows.

\begin{nota}
Let $\hmot$ denote the commutative $\field[\tau]$--algebra
$\field[\tau] [x, y]/ y^2= \tau x$. 
\end{nota}

The algebra $\hmot$ is to be considered as a perturbation of the
algebra $\field [y]$ and the generators $y,x$ are given independent
gradings; this imposes the requirement that the generator $\tau$ has
non-zero grading. In order to encode the non-trivial action on the
underlying ring $\field [\tau]$, one is obliged to work with comodules
over an affine category scheme. The latter is the algebraic object
which represents a small category; it is related to the more familiar
notion of a Hopf algebroid 
(cf  Ravenel \cite[Appendix]{r}) in the same way that 
a bialgebra is related
to a Hopf algebra.

\begin{rem}
An affine category scheme in the category of $\field$--algebras is
given by a pair of commutative $\field$--algebras $(A, \Gamma)$
together with structure morphisms of $\field$--algebras $\eta_L ,
\eta_R \co A \rightrightarrows \Gamma$, $\epsilon \co \Gamma \rightarrow A$,
$\Delta \co \Gamma \rightarrow \Gamma \otimes _ A \Gamma$ which satisfy
a suitable subset of the axioms for a Hopf algebroid.
\end{rem}

In the current situation, the algebra $A$ corresponds to $\field
[\tau]$ and the comodules considered are $\field [\tau]$--modules.

\begin{rem}
The comodules considered in this section need not have underlying
module which is finitely-generated over the algebra $\field
[\tau]$. For this reason, the comodule structures considered are
defined with respect to the completed tensor product, which will
correspond to the underlying grading in applications. The necessary
details are left to be supplied by the reader.
\end{rem}

\begin{prop}
There exists an affine category scheme $(\field [\tau], \cald) $ in
the category of commutative $\field $--algebras, where
 $$
\cald \cong \field[\tau, t, u , \xi_i, \tau_j|i, j \geq 0]/ (u ^2 = t
\xi_0, \tau^2_j = t \tau \xi_{j+1}))
$$
and the structure morphisms are given by $\eta_L \co \tau \mapsto
\tau$, $\eta_R \co \tau \mapsto t \tau $, $\epsilon$ sends $t, u ,\xi_0
$ to $1$ and all other generators to zero. The coproduct $\Delta \co
\cald \rightarrow \cald \otimes_{\field[\tau]} \cald$ is determined by 
\begin{align*}
u &\mapsto u \otimes u &
t &\mapsto t \otimes t \\
\xi_k& \mapsto \sum_{i+j = k}  \xi_j ^{2^i} \otimes \xi_i &
\tau_k& \mapsto \sum_{i+j = k}  \xi_j ^{2^i} \otimes \tau _i +
\tau_k \otimes u .
\end{align*}
Moreover, the $\field[\tau] $--algebra $\hmot$ has the structure of a 
right $(\field[\tau], \cald)$--comodule, with structure morphism 
$\psi \co \hmot \rightarrow \hmot \hctensor_{\field [\tau]} \cald$ which
is the morphism of algebras determined by
$$\tau \mapsto  \tau t \qquad\qquad
y  \mapsto y \otimes u + \sum_{j \geq 0} x ^{2^j} \otimes \tau_j
\qquad\qquad
x \mapsto   \sum_{j \geq 0} x ^{2^j} \otimes \xi_j.$$
\end{prop}

\begin{proof}
The proof is analogous to that of \fullref{prop:calb-construct}. This 
requires the verification that the
only relations imposed by the relation $y^2= \tau x$ are the relations
$u ^2 = t\xi_0$ (which corresponds to the imposed grading upon $t$)
and the relations $\tau_j ^2 = t \tau \xi_{j+1}$.
\end{proof}

\begin{rem}
The above result should be compared with the calculation
of the dual of the algebraic  motivic Steenrod algebra  \cite{v1}
over the coefficient ring $\field [\tau, \rho]$, when $\rho$ is set to
zero. The modification above is that the `grouplike' elements $t, u ,
\xi_0$ are not taken to be $1$. 
\end{rem}

\begin{prop}
There is a unique affine category scheme structure upon 
$$
 (\field [\tau] , \field [\tau, t, u , \xi_0] /  u^2 = t \xi_0)
$$ 
 such that the surjective  morphism of $\field$--algebras $\cald 
\rightarrow
\field [\tau, t, u , \xi_0]/ u^2 = t \xi_0$, defined by sending the 
generators
$\xi_{i+1},  \tau_i$ to zero for $i \geq 1$, induces a morphism of affine
  category schemes 
$$
(\field [\tau], \cald ) \rightarrow (\field [\tau] , 
\field [\tau, t, u , \xi_0] /  u^2 = t \xi_0).
$$ 
\end{prop}

\begin{proof}
Straightforward.
\end{proof}

\begin{rem}
The category of right comodules over the affine category scheme $(\field 
[\tau],\field[\tau,t,u,\xi_0] / u^2 = t \xi_0)$ can be identified with a category of
bigraded $\field[\tau]$--modules, where the generator $\tau$ has
non-zero bidegree.

In the motivic setting, the usual convention is to bigrade so that
 $u$--degree $1$ corresponds to bidegree $(1,1)$, $\xi_0$--degree $1$
 corresponds to bidegree $(2,1)$ and $t$--degree $1$ corresponds to
 bidegree $(0,1)$.
\end{rem}

The definition of the category of graded comodules given in 
\fullref{subsect:graded-comod} generalizes to the current context.

\begin{defn}
Let $\unstmot$ denote the category of graded  right $(\field[\tau] ,
\cald)$--comodules. 
\end{defn}

\begin{prop}
The category $\unstmot$ is a tensor abelian category.
\end{prop}

\subsection[Relating to unstable modules]
{Relating $\unstmot$ to unstable modules} There are
standard base change constructions for affine category schemes;
namely, if $(A, \Gamma)$ is an affine category scheme in the category
of commutative $\field$--algebras and $f \co A\rightarrow B$ is a
morphism of commutative $\field$--algebras, then base change yields an
affine category scheme $(B , B \otimes_A \Gamma \otimes _ A B)$. The
functor on the category of right $A$--modules, $M \mapsto M \otimes _A
B $ extends to a functor from the category of right $(A,
\Gamma)$--comodules to the category of right 
$(B , B \otimes _A \Gamma \otimes _A B)$--comodules, 
which is exact if the functor
$M \mapsto M \otimes _A B$ is exact.
   
The base change constructions apply to the two choices of augmentation
$\epsilon_1, \epsilon_0 \co \field [\tau] \rightrightarrows \field$ given
respectively by $\tau \mapsto 1, \tau \mapsto 0$. To avoid confusion,
the induced bialgebras will be written respectively $\epsilon_1 ^*
\cald $ and $\epsilon_0 ^* \cald$ (noting that an $\field$--affine
category scheme of the form $(\field, \Gamma)$ is a bialgebra). 

\begin{lem}
\label{lem:e1-basechange}
The bialgebra $\epsilon_1 ^* \cald$ is isomorphic to the bialgebra
$\atilde$.
\end{lem}

\begin{proof}
The underlying $\field$--algebra of $\epsilon_1 ^* \cald$ is generated
by the elements $\{t, u , \xi_i , \tau_i \} $ subject to the induced
relations. It is straightforward to verify that the assignment $\tau
\mapsto 1$ implies the identity $t =1$. The result follows readily.
\end{proof}

\begin{defn}
Let $\unstbigrt$ denote the category of graded comodules over
$\epsilon_0 ^* \cald$. 
\end{defn}

\begin{rem}
The category $\unstbigrt$ is not the precise analogue of the category
$\unstbigr$, since the grouplike
element $t$ is present in the bialgebra $\epsilon_0 ^* \cald$ and
there is the relation $u ^2 = t \xi_0$.
\end{rem}

\begin{prop}
The base change functors induce a diagram of functors between tensor
abelian categories
$$
\xymatrix{
\unstmot
\ar[d]_{\epsilon_1 ^*} 
\ar[r]^{\epsilon_0 ^* } 
&
\unstbigrt
\\
\unst
}
$$
in which $\epsilon_1 ^* $ is exact and $\epsilon_0 ^* $ is right
exact. 
\end{prop}

\begin{proof}
The functors are induced by base change, using \fullref{lem:e1-basechange} 
to identify the category of comodules over
$\epsilon_1^* \cald$ with $\unst$. The exactness properties
correspond to the exactness properties of the respective functors $ - 
\otimes
_{\field [\tau] } \field$. 
\end{proof}

\begin{prop}
There is an adjunction 
$$
\epsilon_1 ^* \co \unstmot \rightleftarrows \unst ~{:} \mu. 
$$
Moreover, the underlying $\field [\tau]$--module of $\mu M$, for $M$ an
unstable module, is $\tau$--torsion free. 
\end{prop}

\begin{proof} (Indications) \qua The existence of the right adjoint to
$\epsilon_1 ^*$ is formal, using the existence of a set of
projective generators for the category $\unstmot$. The
$\tau$--torsion statement is a consequence of the surjectivity of the
morphism between projective generators which represents
multiplication by $\tau$.
\end{proof}

\begin{rem}
The existence of $\mu$ and its fundamental properties is implicit in
the work of Yagita \cite{y1} on the motivic cohomology of classifying
spaces of  finite groups.
\end{rem}

\part{Appendix}
\appendix
\section{Comodules} 
\setobjecttype{AppCode}
\label{section-comodules}

Throughout this section, the base ring is taken to be a field, $\field$;
 all algebras are taken to be unital and associative and all coalgebras 
 counital and coassociative. In particular, a bialgebra
will have underlying algebra which is unital and augmented.

\subsection{Duality for comodules and modules}

The following results summarize the elementary properties of duality
in the non-graded case. The category of left modules over an algebra $A$
 is written $_A\module$ and the category of right comodules over a 
coalgebra $C$ is written $\comod _C$; the evident left/right mirror image 
categories are denoted in the obvious way.

\begin{prop}
[Milnor--Moore {\cite[Section 3]{mm}}]
Let $\field$ be a field and let $A$ be an algebra, $C$ a coalgebra and
$B$ a bialgebra over $\field$,  for which the underlying
$\field$--vector spaces are of finite dimension. The following statements
hold:
\begin{enumerate} 
\item
the dual $A^*$ has a natural coalgebra structure;
\item
the dual  $C^*$ has a natural algebra structure;
\item
the canonical morphisms $A \rightarrow A^{** } $ and $C \rightarrow
C^{**}$ are isomorphisms of algebras and coalgebras respectively;
\item
the dual space $B^*$ has the structure of a bialgebra and the
canonical morphism $B \rightarrow B^{**}$ is an isomorphism of
bialgebras;
\item
the bialgebra $B$ has the structure of a Hopf algebra if
and only if the dual $B^*$ has the structure of a Hopf algebra;
\item
there are equivalences of categories: $_A\module  \equiv  \comod _{A^*} $ 
and 
$ \module_A  \equiv  _{A^*} \comod$. 
\end{enumerate}
\end{prop}

The following result is a straightforward application of vector space 
duality.

\begin{prop}
Let $A$ be a finite dimensional $\field$--algebra and let $C$ be a
finite dimensional $\field$--coalgebra. The following statements hold:
\begin{enumerate}
\item
for $M$ a left (respectively right) $A$--module,  the dual vector space 
$M^*$ has  a natural left (resp. right) $A^*$--comodule structure;
\item
for $N$ a left (respectively right) $C$--comodule,  the dual
vector space $N^*$ has a natural left (resp. right) $C^*$--module 
structure.
\end{enumerate}
\end{prop}

The dual $M^*$ of a left module $M$ over an $\field$--algebra $A$  has the
structure of a right $A$--module; under a finiteness hypothesis, there
is an analogue for coalgebras.

\begin{prop}
Let $C$ be an $\field$--coalgebra of finite dimension and let $N$ be a
left $C$--comodule, then the dual vector space  $N^*$ has the structure
of a right $C$--comodule. 
\end{prop}

\begin{proof}
The right $C$--comodule structure is adjoint to the left
$C^*$--module structure which is given by the dual $C^* \otimes N^*
\rightarrow N^*$ to the comodule structure morphism.
\end{proof}

\subsection{Corestriction and cotensor products}
\setobjecttype{AppCode}
\label{sect:cotensor}

The following standard result corresponds to the fact that
the category of modules over a ring is abelian. 

\begin{prop}
The category of right (respectively left) comodules over a
$\field$--coalgebra  is an abelian category.
\end{prop}

The cotensor product of two comodules over a coalgebra is formally dual to 
the
definition of the tensor product of two modules over an algebra.

\begin{defn}
Let $C$ be a $\field$--coalgebra and let $M$ be a right $C$--comodule
and let $N$ be a left $C$--comodule;  the cotensor product $M\Box
_C N$ is the kernel of the morphism $\psi_M \otimes 1 - 1 \otimes
\psi_N \co M\otimes N \rightarrow M\otimes C\otimes N$, where $\psi_M,
\psi_N$ denote the respective structure morphisms of $M, N$. 
\end{defn}

\begin{lem}
\label{lem-cotensor-C}
For $M \in \comod_C$  a right $C$--comodule,  there is a natural
isomorphism of $\field$--vector spaces $M \cong M \Box _C C$.
\end{lem}

\begin{proof}
This fundamental result is a consequence of the counital axiom for
the coalgebra $C$.
\end{proof}

The following result is standard.

\begin{prop}
 Let $C$ be an $\field$--coalgebra and let $N$ be a left
$C$--comodule, then the cotensor product $- \Box_C N \co \comod_C
\rightarrow \vs$ is a left exact functor.
\end{prop}

\begin{exam}
Let $C$ denote the underlying coalgebra of the $\field$--Hopf algebra
$\field[u,u^{-1}]$ and let $\field (n)$ denote the left
$C$--comodule $\field$ with structure morphism $\field \rightarrow
C \otimes \field $, $1 \mapsto  x^n\otimes 1$, for some integer $n$. The
functor $- \Box _C \field (n)$ is an exact functor from the category
of right $C$--comodules to the category of vector spaces. This functor
corresponds to the projection of a graded vector space onto the
component of degree $n$.
\end{exam}

The following result has an obvious generalization to the graded
finite-type case, when a suitable connectivity hypothesis is imposed on the
graded modules (cf \cite[Proposition 3.2]{mm}, where all graded
modules are taken to be connective in a suitable sense).

\begin{lem}
Let $C$ be an $\field$--coalgebra of finite dimension, let $M$ be a
right $C$--comodule and let $N$ be a left $C$--comodule, such that $M,N$ 
are of finite dimension. There is
an isomorphism of vector spaces 
$
(M\Box_C N) ^* \cong M^* \otimes _{C^*} N^* $.
\end{lem}

\begin{defn}
For  $f \co D \rightarrow C$  a morphism of $\field$--coalgebras, the
corestriction functor $f_* \co \comod_D \rightarrow \comod_C$ is the
functor which sends $M \in \comod_D$ to the right $C$--comodule $M \in 
\comod_C$, with  structure morphism the
composite
$$M \rightarrow M\otimes D \stackrel{1\otimes f}{\longrightarrow}
  M \otimes C.$$ 
The corestriction functor   $f_* \co \  _D\comod
\rightarrow \ _C\comod$ is defined analogously. 
\end{defn}

\begin{prop}
Let $f: D \rightarrow C$ be a morphism of $\field$--coalgebras, then
the corestriction functor  $f_* \co \comod_D \rightarrow \comod_C$
admits a right adjoint $f^! \co \comod_C \rightarrow \comod_D$ given by
$N \in \comod _C \mapsto N \Box_C D$, where the $D$--comodule structure
on $N\Box_C D$ is induced by the coproduct $D \rightarrow D \otimes
D$. 
\end{prop}

\begin{proof}
The adjunction morphisms are induced by the following constructions. For 
$M$ a right
$D$--comodule,  there is a canonical morphism $M \rightarrow (f_*
M) \Box _C D $ which is induced by the structure morphism $M
\rightarrow M \otimes D$. For $N$ a right $C$--comodule,  the
counit of the adjunction is the morphism $N\Box_C D \rightarrow N\Box
_C C \cong C$ which is induced by $f$, where the isomorphism is
provided by \fullref{lem-cotensor-C}.
\end{proof}

The corestriction functor admits a left adjoint when the coalgebras
satisfy suitable finite-type hypotheses. The ungraded version of the
result is the following:

\begin{prop}
Let $f: D \rightarrow C$ be a morphism of $\field$--coalgebras, where $D,
C$ are of finite dimension. The corestriction functor $f_* \co
\comod_D \rightarrow \comod_C$ admits a left adjoint $f^* \co \comod_C
\rightarrow \comod_D$ given by $N \mapsto D^* \otimes _{C^*} N$,
where $N$ is regarded as a left $C^*$--module and the right
$D$--comodule structure is adjoint to the extended left $D^*$--module
structure.
\end{prop}

\begin{rem} 
Suppose that the right $C$--comodule $N$ is of finite dimension, where
$f \co D\rightarrow C$ satisfies the hypotheses of the Proposition, then
there is an isomorphism of right $D$--comodules $(D \Box_C N^*) ^*
\cong D^* \otimes _{C^*} N$. This result generalizes to the context
of graded comodules if the graded objects $N,C ,D$ are all of finite
type and the dual algebras $D^* , C^*$ are connective (trivial in
sufficiently negative dimensions).
\end{rem}

\bibliographystyle{gtart}
\bibliography{link}
\end{document}